\documentclass[titlepage,12pt]{article} 
\usepackage{pstricks,pst-plot} 
\usepackage[title]{appendix}
\usepackage{amssymb,amsthm,amsmath} 
\usepackage[a4paper]{geometry}
\usepackage{enumitem}
\usepackage[utf8]{inputenc}
\usepackage{hyperref}

\date{}

\newcommand{\n}{\mathbb{N}}
\newcommand{\z}{\mathbb{Z}}

\newcommand{\re}{\mathbb{R}}

\newcommand{\ep}{\varepsilon}
\renewcommand{\qed}{{\penalty 10000\mbox{$\quad\Box$}}}

\newcommand{\li}{\lambda_{n}}
\newcommand{\liu}{\lambda_{n-1}}
\newcommand{\ul}{u_{\lambda}}
\newcommand{\G}{\mathcal{G}}
\newcommand{\cka}{C^{k,\alpha}}
\newcommand{\gl}{\gamma_{\lambda}}
\newcommand{\eg}{E_{\gamma}}
\newcommand{\egl}{E_{\gamma_{\lambda}}}
\newcommand{\cep}{c_{\ep}}
\newcommand{\ceph}{\widehat{c}_{\ep}}

\newtheorem{thm}{Theorem}[section]
\newtheorem{thmbibl}{Theorem}

\newtheorem{rmk}[thm]{Remark}
\newtheorem{prop}[thm]{Proposition}
\newtheorem{defn}[thm]{Definition}

\newtheorem{lemma}[thm]{Lemma}
\newtheorem*{open}{Open problem}

\title{Time-dependent propagation speed vs strong damping for degenerate linear hyperbolic equations}

\author{Marina Ghisi\vspace{1ex}\\ 
{\normalsize Universit\`a degli Studi di Pisa} \\
{\normalsize Dipartimento di Matematica}\\ 
{\normalsize PISA (Italy)}\\
{\normalsize e-mail: \texttt{marina.ghisi@unipi.it}}
\and
Massimo Gobbino\vspace{1ex}\\ 
{\normalsize Universit\`a degli Studi di Pisa} \\
{\normalsize Dipartimento di Ingegneria Civile e Industriale}\\ 
{\normalsize PISA (Italy)}\\  
{\normalsize e-mail: \texttt{massimo.gobbino@unipi.it}}
}

\begin{document}
\maketitle
\begin{abstract}

We consider a degenerate abstract wave equation with a time-dependent propagation speed. We investigate the influence of a strong dissipation, namely a friction term that depends on a power of the elastic operator.

We discover a threshold effect. If the propagation speed is regular enough, then the damping prevails, and therefore the initial value problem is well-posed in Sobolev spaces. Solutions also exhibit a regularizing effect analogous to parabolic problems. As expected, the stronger is the damping, the lower is the required regularity.

On the contrary, if the propagation speed is not regular enough, there are examples where the damping is ineffective, and the dissipative equation behaves as the non-dissipative one.

\vspace{6ex}

\noindent{\bf Mathematics Subject Classification 2010 (MSC2010):} 35L20, 35L80, 35L90.


\vspace{6ex}


\end{abstract}

 
\section{Introduction}

Let $H$ be a separable real Hilbert space.  For every $x$ and $y$ in $H$, $|x|$ denotes the norm of $x$, and $\langle x,y\rangle$ denotes the scalar product of $x$ and $y$.  Let $A$ be a self-adjoint linear operator on $H$ with dense domain $D(A)$.  We assume that $A$ is nonnegative, namely $\langle Ax,x\rangle\geq 0$ for every $x\in D(A)$, so that for every $\alpha\geq 0$ the power $A^{\alpha}x$ is defined provided that $x$ lies in a suitable domain $D(A^{\alpha})$.

We consider the second order linear evolution equation
\begin{equation}
	u''(t)+2\delta A^{\sigma}u'(t)+c(t)Au(t)=0
	\label{pbm:eqn}
\end{equation}
in some interval $[0,T]$, with initial data
\begin{equation}
	u(0)=u_{0},
	\hspace{3em}
	u'(0)=u_{1}.
	\label{pbm:data}
\end{equation}

We refer to~\cite{gg:dgcs-strong} for the history of the problem and a short survey of some related literature (see also~\cite{Bui-Reissig,DAbbicco-Ebert,Ebert-Reissig,Hirosawa-Bui} and the references quoted therein for analogous models with competition between damping and time-dependent propagation speed). Here we just recall the main results that are more relevant to our presentation. 

The non-dissipative equation ($\delta=0$) was considered in the seminal paper~\cite{dgcs} under the strict hyperbolicity assumption
\begin{equation}
0<\mu_{1}\leq c(t)\leq\mu_{2}
	\qquad
	\forall t\in[0,T],
\label{hp:sh}
\end{equation}
and then in~\cite{cjs} under the degenerate hyperbolicity assumption
\begin{equation}
	0\leq c(t)\leq\mu
	\qquad
	\forall t\in[0,T].
	\label{hp:wh}
\end{equation}

The general philosophy is that higher space-regularity of initial data compensates lower time-regularity of $c(t)$. The result is that problem (\ref{pbm:eqn})--(\ref{pbm:data}) is well-posed in suitable Gevrey spaces, whose order depends on the regularity class of $c(t)$, and on the strict/degenerate hyperbolicity condition. For less regular data strange pathologies may occur, in the sense that for suitable coefficients there do exist ``solutions'' which lie in Gevrey spaces (of course not as good as those that guarantee well-posedness) at time $t=0$, but which are not even distributions when $t>0$. We refer to section~\ref{sec:thmbibl} for a survey of the statements concerning the degenerate non-dissipative case.

The dissipative equation ($\delta>0$) with constant positive propagation speed was considered in full generality in~\cite{ggh:strong-damping}. If we limit ourselves to the range $\sigma\in(0,1/2]$, the result is that in this special autonomous case problem (\ref{pbm:eqn})--(\ref{pbm:data}) is well-posed in the classic energy space $D(A^{1/2})\times H$, and solutions exhibit a regularizing effect for positive times, in the sense that they lie in Gevrey spaces of order $(2\sigma)^{-1}$.

The dissipative case with time-dependent propagation speed is more complex, because there is some sort of competition between the damping and the potential low-regularity of $c(t)$. This competition was investigated for the first time in~\cite{gg:dgcs-strong}, leading to the following results. 

\begin{itemize}

\item  When $\sigma>1/2$ the damping always prevails, and problem (\ref{pbm:eqn})--(\ref{pbm:data}) is well-posed in $D(A^{1/2})\times H$ (but also different choices are possible) provided that $c(t)$ is measurable and satisfies the degenerate hyperbolicity condition (\ref{hp:wh}).

\item  When $\sigma\in[0,1/2]$ the competition is tighter. If $c(t)$ is $\alpha$-H\"older continuous and satisfies the strict hyperbolicity condition (\ref{hp:sh}), then problem (\ref{pbm:eqn})--(\ref{pbm:data}) is well-posed in $D(A^{1/2})\times H$ provided that $2\sigma>1-\alpha$. Otherwise, the equation behaves as the non-dissipative one, meaning well-posedness in the appropriate Gevrey classes, and potential pathologies for less regular data.

\end{itemize}

In this paper, which is intended as a continuation of~\cite{gg:dgcs-strong}, we consider the case where $\sigma\in[0,1/2]$, and the coefficient $c(t)$ is a function of class $C^{k,\alpha}$ satisfying the degenerate hyperbolicity condition (\ref{hp:wh}). Again we discover a threshold effect. 
\begin{itemize}

\item When $(2+k+\alpha)\sigma>1$, we show in Theorem~\ref{thm:main} that equation (\ref{pbm:eqn}) behaves as the one with constant positive propagation speed, meaning well-posedness in Sobolev spaces, and regularizing effect to Gevrey classes of order $(2\sigma)^{-1}$ for positive times.

\item  When $(2+k+\alpha)\sigma<1$, we show in Theorem~\ref{thm:dgcs} that equation (\ref{pbm:eqn}) can exhibit the same pathologies of the non-dissipative case.

\end{itemize}

From the technical point of view, the spectral theory reduces the problem to estimating the growth of solutions to the family of ordinary differential equations 
\begin{equation}
	\ul''(t)+2\delta\lambda^{2\sigma}\ul'(t)+
	\lambda^{2}c(t)\ul(t)=0,
	\label{pbm:main-ode}
\end{equation}
with initial data
\begin{equation}
	\ul(0)=u_{0,\lambda},
	\hspace{3em}
	\ul'(0)=u_{1,\lambda}.
	\label{eqn:ODE-data}
\end{equation}

To this end, we introduce ``approximated hyperbolic energies'' of the form
\begin{equation}
|\ul'(t)|^{2}+\delta^{2}\lambda^{4\sigma}|\ul(t)|^{2}+\delta\lambda^{2\sigma}\ul(t)\ul'(t)+\gl(t)\lambda^{2}|\ul(t)|^{2},
\nonumber
\end{equation}
where $\gl(t)$ is a suitable smooth approximation of $c(t)$ to be chosen in a $\lambda$-dependent way. This technique dates back to~\cite{dgcs,cjs}, but here we need to design $\gl(t)$ in a completely different way in order to take advantage of the strong damping. For this reason, Lemma~\ref{lemma:gamma} is the technical core of the proof of Theorem~\ref{thm:main}.

As for counterexamples, again we follow the strategy devised in~\cite{dgcs,cjs}, but again we have to change the ingredients from the very beginning because of the dissipation. 
 
This paper is organized as follows.  In section~\ref{sec:notation} we introduce the functional setting and we recall the classic existence results from~\cite{cjs}.  In section~\ref{sec:main} we state our main results.  In section~\ref{sec:heuristics} we provide a heuristic description of the competition between oscillations of $c(t)$ and strong damping.  In section~\ref{sec:proofs} we prove our existence and regularity results.  In section~\ref{sec:counterexamples} we present our examples of pathological solutions.


\setcounter{equation}{0}
\section{Notation and previous results}\label{sec:notation}

\subsection{Functional spaces}

Let $H$ be a separable Hilbert space.  Let us assume that $H$ admits a
countable complete orthonormal system $\{e_{n}\}_{n\in\n}$ made by
eigenvectors of $A$.  We denote the corresponding eigenvalues by
$\lambda_{n}^{2}$ (with the agreement that $\lambda_{n}\geq 0$), so that
$Ae_{n}=\lambda_{n}^{2}e_{n}$ for every $n\in\n$.  In this case every
$u\in H$ can be written in a unique way in the form
$u=\sum_{n=0}^{\infty}u_{n}e_{n}$, where $u_{n}=\langle
u,e_{n}\rangle$ are the Fourier components of $u$.  In other words,
the Hilbert space $H$ can be identified with the set of sequences
$\{u_{n}\}$ of real numbers such that
$\sum_{n=0}^{\infty}u_{n}^{2}<+\infty$.

We stress that this is just a simplifying assumption, with
substantially no loss of generality.  Indeed, according to the
spectral theorem in its general form (see for example Theorem~VIII.4
in~\cite{reed}), one can always identify $H$ with $L^{2}(M,\mu)$ for a
suitable measure space $(M,\mu)$, in such a way that under this
identification the operator $A$ acts as a multiplication operator by
some measurable function $\lambda^{2}(\xi)$.  All definitions and
statements in the sequel, with the exception of the counterexamples
of Theorem~\ref{thm:dgcs}, can be easily extended to the
general setting just by replacing the sequence $\{\lambda_{n}^{2}\}$ with the
function $\lambda^{2}(\xi)$, and the sequence $\{u_{n}\}$ of Fourier
components of $u$ with the element $\widehat{u}(\xi)$ of
$L^{2}(M,\mu)$ corresponding to $u$ under the identification of $H$
with $L^{2}(M,\mu)$.

The usual functional spaces can be characterized in terms of Fourier
components as follows.

\begin{defn}
	\begin{em}
		Let $u$ be a sequence $\{u_{n}\}$ of real numbers.
		\begin{itemize}
			\item  \emph{Sobolev spaces}. For every $\alpha\geq 0$ we say that $u\in D(A^{\alpha})$ if 
			\begin{equation}
				\|u\|_{D(A^{\alpha})}^{2}:=
				\sum_{n=0}^{\infty}(1+\lambda_{n})^{4\alpha}u_{n}^{2}<+\infty.
				\label{defn:sobolev}
			\end{equation}
		
			\item \emph{Distributions}.  We say that $u\in
			D(A^{-\alpha})$ for some $\alpha\geq 0$ if 
			\begin{equation}
				\|u\|_{D(A^{-\alpha})}^{2}:=
				\sum_{n=0}^{\infty}(1+\lambda_{n})^{-4\alpha}u_{n}^{2}<+\infty.
				\label{defn:distributions}
			\end{equation}
		
\item \emph{Gevrey spaces}.  Let $s>0$, $r>0$ and $\alpha$ be real numbers. We say that $u\in\G_{s,r,\alpha}(A)$ if 
\begin{equation}
	\|u\|^{2}_{s,r,\alpha}:=
	\sum_{n=0}^{\infty}(1+\lambda_{n})^{4\alpha}u_{n}^{2} \exp\left(2r\lambda_{n}^{1/s}\right)<+\infty.
	\label{defn:gevrey}
\end{equation}
		
\item \emph{Gevrey ultradistributions}.  Let $S>0$, $R>0$ and $\alpha$ be real numbers. We say that $u\in\G_{-S,R,\alpha}(A)$ if 
\begin{equation}
	\|u\|^{2}_{-S,R,\alpha}:=
	\sum_{n=0}^{\infty}(1+\lambda_{n})^{4\alpha}u_{n}^{2} \exp\left(-2R\lambda_{n}^{1/S}\right)<+\infty.
	\label{defn:hyperfunctions}
\end{equation}

\end{itemize}

	\end{em}
\end{defn}

The quantities defined in (\ref{defn:sobolev}) through (\ref{defn:hyperfunctions}) are actually norms inducing a Hilbert space structure on the corresponding spaces.  The standard inclusions
$$\G_{s,r,\alpha}(A)\subseteq D(A^{\beta})\subseteq H\subseteq D(A^{-\beta})\subseteq\G_{-S,R,-\alpha}(A)$$
hold true for every positive value of $\alpha$, $\beta$, $r$, $s$, $R$, and $S$.  All inclusions are strict if the sequence $\lambda_{n}$ is unbounded.

We observe that $\G_{s,r,\alpha}(A)$ is actually a so-called \emph{scale of Hilbert spaces} with respect to the parameter $r$, with larger values of $r$ corresponding to smaller spaces.  Analogously, $\G_{-S,R,\alpha}(A)$ is a scale of Hilbert spaces with respect to the parameter $R$, but with larger values of $R$ corresponding to larger spaces.

\subsection{Damping-independent results}\label{sec:thmbibl}

In this subsection we recall the classical results concerning existence, uniqueness, and regularity for solutions to problem (\ref{pbm:eqn})--(\ref{pbm:data}) under the sole assumptions that $\delta\geq 0$ and $c(t)$ satisfies the degenerate hyperbolicity assumption. For the sake of consistency, we rephrase the results in our functional setting. In the quoted references only the case $\delta=0$ is considered, but the same techniques work also when $\delta> 0$ because all extra terms have the ``right sign''. 

The first result concerns existence and uniqueness of a very weak solution for a very huge class of initial data, with minimal assumptions on $c(t)$ (no hyperbolicity is required).

\begin{thmbibl}[see~{\cite[Theorem~1]{dgcs}}]\label{thmbibl:dh}

Let us consider problem (\ref{pbm:eqn})--(\ref{pbm:data}) under the following assumptions:
\begin{itemize}

\item $A$ is a self-adjoint nonnegative operator on a separable Hilbert space $H$,
	
\item $c\in L^{1}((0,T))$ (without sign conditions) for some $T>0$,
	
\item $\sigma\geq 0$ and $\delta\geq 0$ are two real numbers,
	
\item there exists $R_{0}>0$ such that initial conditions satisfy 
$$(u_{0},u_{1})\in\G_{-1,R_{0},1/2}(A)\times\G_{-1,R_{0},0}(A).$$ 

\end{itemize}
	
Then there exists a nondecreasing function $R:[0,T]\to(0,+\infty)$, with $R(0)=R_{0}$, such that problem (\ref{pbm:eqn})--(\ref{pbm:data}) admits a unique solution
\begin{equation}
	u\in C^{0}\left([0,T];\G_{-1,R(t),1/2}(A)\right)\cap C^{1}\left([0,T];\G_{-1,R(t),0}(A)\right).		
	\label{th:(DGCS)-psi}
\end{equation}

\end{thmbibl}

Condition (\ref{th:(DGCS)-psi}), with the range space increasing with time,
simply means that 
$$u\in
C^{0}\left([0,\tau];\G_{-1,R(\tau),1/2}(A)\right)\cap
C^{1}\left([0,\tau];\G_{-1,R(\tau),0}(A)\right)
\quad\quad
\forall\tau\in(0,T].$$
		
This amounts to say that scales of Hilbert spaces, rather than fixed Hilbert spaces, are the natural setting for this problem.

In the second result we assume degenerate hyperbolicity and more time-regularity of the coefficient $c(t)$, and we obtain well-posedness in a suitable smaller class of Gevrey ultradistributions.

\begin{thmbibl}[see~{\cite[Theorem~1 and Remark~4]{cjs}}]\label{thmbibl:existence}

Let us consider problem (\ref{pbm:eqn})--(\ref{pbm:data}) under the following assumptions:
\begin{itemize}

\item $A$ is a self-adjoint nonnegative operator on a separable Hilbert space $H$,
	
\item  there exists $k\in\n$ and $\alpha\in(0,1]$ such that $c\in\cka([0,T])$,
	
\item  $c(t)$ satisfies the degenerate hyperbolicity assumption (\ref{hp:wh}), 

\item $\sigma\geq 0$ and $\delta\geq 0$ are two real numbers,
	
\item initial conditions satisfy
$$(u_{0},u_{1})\in\G_{-S,R_{0},1/2}(A)\times\G_{-S,R_{0},0}(A)$$ 
for some real numbers $R_{0}>0$ and $S>0$ such that
\begin{equation}
	S<1+\frac{k+\alpha}{2}.
	\label{hp:Ska}
\end{equation}

\end{itemize}
	
Then the unique solution $u(t)$ to the problem provided by Theorem~\ref{thmbibl:dh} satisfies the further regularity
$$u\in C^{0}\left([0,T],\G_{-S,R_{0}+\ep,1/2}(A)\right)\cap C^{1}\left([0,T],\G_{-S,R_{0}+\ep,0}(A)\strut\right)
\qquad
\forall\ep>0.$$
	
\end{thmbibl}

The third result concerns existence of regular solutions. The assumptions on $c(t)$ are the same as in Theorem~\ref{thmbibl:existence}, but initial data are significantly more regular (Gevrey spaces instead of Gevrey ultradistributions).

\begin{thmbibl}[see~{\cite[Theorem~1]{cjs}}]\label{thmbibl:regularity}

Let us consider problem (\ref{pbm:eqn})--(\ref{pbm:data}) under the following assumptions:
\begin{itemize}

\item $A$ is a self-adjoint nonnegative operator on a separable Hilbert space $H$,
	
\item  there exists $k\in\n$ and $\alpha\in(0,1]$ such that $c\in\cka([0,T])$,
	
\item  $c(t)$ satisfies the degenerate hyperbolicity assumption (\ref{hp:wh}), 

\item $\sigma\geq 0$ and $\delta\geq 0$ are two real numbers,
	
\item initial conditions satisfy
$$(u_{0},u_{1})\in\G_{s,r_{0},1/2}(A)\times\G_{s,r_{0},0}(A)$$ 
for some real numbers $r_{0}>0$ and $s>0$ such that 
\begin{equation}
	s<1+\frac{k+\alpha}{2}.
	\label{hp:ska}
\end{equation}

\end{itemize}
	
Then the unique solution $u(t)$ to the problem provided by Theorem~\ref{thmbibl:dh} satisfies the further regularity
\begin{equation}
	u\in C^{0}\left([0,T],\G_{s,r_{0}-\ep,1/2}(A)\right)\cap C^{1}\left([0,T],\G_{s,r_{0}-\ep,0}(A)\strut\right)
	\qquad
	\forall\ep\in(0,r_{0}).
\nonumber
\end{equation}
	
\end{thmbibl}

\begin{rmk}\label{rmk:cjs}
\begin{em}

The counterexample presented in~\cite[Theorem~2]{cjs} clarifies that there is essentially no well-posedness result in between the Gevrey spaces of Theorem~\ref{thmbibl:regularity} and the Gevrey ultradistributions of Theorem~\ref{thmbibl:existence}, and that conditions (\ref{hp:Ska}) and (\ref{hp:ska}) are optimal. 

More precisely, there exists a nonnegative coefficient $c(t)$ of class $C^{k,\alpha}$ for which (\ref{pbm:eqn}) admits a solution that is Gevrey regular at time $t=0$ (just a little bit less regular than required by Theorem~\ref{thmbibl:regularity}), but then exhibits a severe derivative loss, meaning that for all positive times this solution is just a hyperdistribution as in Theorem~\ref{thmbibl:existence}, and nothing more. 

\end{em}
\end{rmk}


\subsection{Glaeser type inequalities}

A classical result states that the power $1/(k+\alpha)$ of a nonnegative function of class $C^{k,\alpha}$ is absolutely continuous, and actually Lipschitz continuous when $k=1$. In this paper we need this result in the following form.

\begin{thmbibl}[Glaeser type inequalities]\label{thmbibl:glaeser}

Let $T$ be a positive real number, let $k$ be a positive integer, let $\alpha\in(0,1]$ be a real number, and let $c:[0,T]\to[0,+\infty)$ be a nonnegative function of class $C^{k,\alpha}$.

Then the following estimates hold true.
\begin{itemize}
  \item \emph{(Case $k=1$)} There exists a constant $K$ such that
  \begin{equation}
  |c'(t)|\leq K[c(t)]^{1-1/(1+\alpha)}
  \qquad
  \forall t\in[0,T].
  \label{th:glaeser-1}
  \end{equation}
  
  \item \emph{(Case $k\geq 2$)} There exists a function $\varphi:[0,T]\to[0,+\infty)$, with $\varphi\in L^{1}((0,T))$, such that
  \begin{equation}
  |c'(t)|\leq \varphi(t)[c(t)]^{1-1/(k+\alpha)}
  \qquad
  \forall t\in[0,T].
  \label{th:glaeser-2}
  \end{equation}

\end{itemize}

\end{thmbibl}

For a proof of Theorem~\ref{thmbibl:glaeser} we refer to~\cite[Lemma~1]{cjs}, or to the more recent paper~\cite{gg:roots} where the result has been improved by showing further $L^{p}$ summability of $\varphi(t)$. 


\setcounter{equation}{0}
\section{Main results}\label{sec:main}

Let us set
\begin{equation}
C(t):=\int_{0}^{t}c(s)\,ds
\qquad
\forall t\in[0,T].
\label{defn:C}
\end{equation}

Our main existence and regularity result concerns the regime where the damping dominates the time-dependent coefficient.

\begin{thm}[Sobolev and Gevrey regularity]\label{thm:main}

Let us consider problem (\ref{pbm:eqn})--(\ref{pbm:data}) under the following assumptions:
\begin{itemize}
		
\item $A$ is a self-adjoint nonnegative operator on a separable Hilbert space $H$,
	
\item  there exists $k\in\n$ and $\alpha\in(0,1]$ such that $c\in\cka([0,T])$,
	
\item  $c(t)$ satisfies the degenerate hyperbolicity assumption (\ref{hp:wh}), with in addition
\begin{equation}
c(0)=0,
\label{hp:c0=0}
\end{equation} 
and its antiderivative (\ref{defn:C}) satisfies
\begin{equation}
C(t)>0
\qquad
\forall t\in(0,T],
\label{hp:C>0}
\end{equation}

\item $\delta$ is a positive real number, and $\sigma$ is a real number such that
\begin{equation}
	\frac{1}{2+k+\alpha}<\sigma\leq\frac{1}{2}.	
	\label{hp:sigma-ka}
\end{equation}
				
\item $(u_{0},u_{1})\in D(A^{\sigma})\times H$.

\end{itemize}

	Then the unique solution $u(t)$ to the problem provided by Theorem~\ref{thmbibl:dh} has the following regularity properties.

\begin{enumerate}
\renewcommand{\labelenumi}{(\arabic{enumi})}

\item \emph{(Sobolev regularity for positive times)} It turns out that
\begin{equation}
u\in C^{0}\left((0,T],D(A^{\sigma})\right)\cap C^{1}\left((0,T],H\right).
\label{th:main-sob-open}
\end{equation}

\item \emph{(Gevrey regularity for positive times)} There exists $r>0$ such that
\begin{equation}
u\in C^{0}\left((0,T],\G_{(2\sigma)^{-1},rC(t),\sigma}(A)\right)\cap C^{1}\left((0,T],\G_{(2\sigma)^{-1},rC(t),0}(A)\right).
\label{th:main-gevrey}
\end{equation}

\item \emph{(Continuity in Sobolev spaces up to $t=0$)} If in addition $k\in\{0,1\}$, then it turns out that
\begin{equation}
u\in C^{0}\left([0,T],D(A^{\sigma})\right)\cap C^{1}\left([0,T],H\right).
\label{th:main-sob-closed}
\end{equation}

\end{enumerate}	

\end{thm}

The proof of Theorem~\ref{thm:main} provides also estimates for high frequency components of solutions. We refer to Remark~\ref{rmk:estimates} for further details.

We observe that we assumed that $c(0)=0$, and that $c(t)$ does not vanish identically in a right neighborhood of $t=0$. In the following two remarks we show that there is no loss of generality in these assumptions.

\begin{rmk}
\begin{em}

Let us assume that $c(0)>0$. Due to the continuity of $c(t)$, there exists $T_{1}\in(0,T]$ such that $c(t)$ is bounded from below by a positive constant in $[0,T_{1}]$. 

Therefore, in this subinterval we can apply the theory for the strictly hyperbolic case, which provides well-posedness in Sobolev spaces (see~\cite[Theorem~3.2 and Remark~3.5]{gg:dgcs-strong}) and regularizing effect up to Gevrey spaces of order $(2\sigma)^{-1}$ (see~\cite[Theorem~3.9]{gg:dgcs-strong}) provided that $c\in C^{0,\beta}$ for some $\beta>1-2\sigma$. In order to check this assumption, we observe that it is satisfied with $\beta:=\alpha$ if (\ref{hp:sigma-ka}) holds true with $k=0$, and with $\beta$ close enough to~1 if (\ref{hp:sigma-ka}) holds true with $k\geq 1$ (in which case $c(t)$ is at least of class $C^{1}$, and hence $\beta$-H\"older continuous for every $\beta\in(0,1)$).

Thus from the theory for strictly hyperbolic equations we deduce that, if initial data are in $D(A^{1/2})\times H$, there exists $r_{1}>0$ such that
$$(u(T_{1}),u'(T_{1}))\in\G_{s,r_{1},1/2}(A)\times\G_{s,r_{1},0}(A)
\qquad\mbox{with $s=(2\sigma)^{-1}$.}$$

This value of $s$ satisfies (\ref{hp:ska}) because of (\ref{hp:sigma-ka}), and hence we can apply Theorem~\ref{thmbibl:regularity} in the interval $[T_{1},T]$, where $c(t)$ is allowed to vanish. As a consequence, a quite regular solution exists on the whole interval $[0,T]$ for all initial data in $D(A^{1/2})\times H$.

In other words, the case $c(0)>0$ can be dealt with relying only on Theorem~\ref{thmbibl:regularity} and on the theory for strictly hyperbolic equations.

\end{em}
\end{rmk}

\begin{rmk}\label{rmk:c==0}
\begin{em}

Let us assume that $c(t)$ vanishes identically in a left neighborhood of the origin, and let us set
$$T_{1}:=\sup\{t\in[0,T]:c(\tau)=0\quad\forall\tau\in[0,t]\}.$$

In the interval $[0,T_{1}]$ equation (\ref{pbm:main-ode}) reduces to
$$\ul''(t)+2\delta\lambda^{2\sigma}\ul'(t)=0,$$
whose solution is
\begin{equation}
\ul(t)=u_{0,\lambda}-\frac{\exp(-2\delta\lambda^{2\sigma}t)-1}{2\delta\lambda^{2\sigma}}\cdot u_{1,\lambda}
\qquad
\forall t\in[0,T_{1}].
\label{sol:c=0}
\end{equation}

If $(u_{0},u_{1})\in D(A^{\sigma})\times H$, this formula tells us that
$$u\in C^{0}\left([0,T_{1}],D(A^{\sigma})\right)\cap C^{1}\left([0,T_{1}],H\right),$$
and in particular $(u(T_{1}),u'(T_{1}))\in D(A^{\sigma})\times H$. Therefore, we are in a position to continue the solution by applying Theorem~\ref{thm:main} in the interval $[T_{1},T]$.

\end{em}
\end{rmk}

\begin{rmk}
\begin{em}

The calculation shown in Remark~\ref{rmk:c==0} clarifies that $D(A^{\sigma})\times H$ (or any product of spaces with ``gap $\sigma$'') is the appropriate phase space for degenerate equations. In some sense, this space is chosen by the equation itself. Indeed, formula (\ref{sol:c=0}) with $u(0)=0$ and $u'(0)=u_{1}$ shows that solutions can undergo an immediate jump
$$D(A^{\infty})\times H	\rightsquigarrow D(A^{\sigma})\times H.$$

As expected, this sort of derivative loss is bigger when $\sigma$ is smaller. 

\end{em}
\end{rmk}

\begin{rmk}
\begin{em}

In Theorem~\ref{thm:main} we prove the continuity of the solution in $D(A^{\sigma})\times H$ up to $t=0$ only in the case $k\leq 1$. When $k\geq 2$, we obtain Sobolev and Gevrey regularity for positive times, but as far as we know the solution might assume initial data only in the weak hyperdistributional sense of Theorem~\ref{thmbibl:existence}, even if these initial data are again in $D(A^{\sigma})\times H$.

From the technical point of view, this depends on the fact that for $k\geq 2$ the function $\varphi$ that appears in the Glaeser type inequalities of Theorem~\ref{thmbibl:glaeser} could be unbounded. We refer to Remark~\ref{rmk:ka} for further details. 

On the other hand, we have no counterexamples to the continuity up to $t=0$, which motivates the following question.

\end{em}
\end{rmk}

\begin{open}

Let us consider problem (\ref{pbm:eqn})--(\ref{pbm:data}) under the same assumptions of Theorem~\ref{thm:main}. Can we conclude that (\ref{th:main-sob-closed}) holds true even in the case $k\geq 2$?

\end{open}

Our second result is the counterpart of Theorem~\ref{thm:main}, and concerns the regime where the damping is ineffective. In this regime Theorem~\ref{thmbibl:regularity} still provides existence of a regular solution for initial data in suitable Gevrey classes. Here we show that for less regular data a severe derivative loss is possible.

\begin{thm}[Instantaneous severe derivative loss]\label{thm:dgcs}

Let $A$ be a linear operator on a Hilbert space $H$. Let us assume that there exists a countable (not necessarily complete) orthonormal system $\{e_{n}\}$ in $H$, and an unbounded sequence $\{\lambda_{n}\}$ of positive real numbers such that $Ae_{n}=\lambda_{n}^{2}e_{n}$ for every $n\in\n$. 
	
Let  $\delta>0$ and $\alpha\in(0,1]$ be real numbers, let $k\in\n$ be a nonnegative integer, and let $\sigma$ be a real number such that
\begin{equation}
0\leq\sigma<\frac{1}{2+k+\alpha}.
\label{hp:sigma-ka-cjs}
\end{equation}

Then there exist a function $c:\re\to [0,+\infty)$, and a solution $u(t)$ to equation (\ref{pbm:eqn}) in $[0,+\infty)$ satisfying the following three properties.
\begin{enumerate}
\renewcommand{\labelenumi}{(\arabic{enumi})}

\item (Regularity of the coefficient) The coefficient $c(t)$ satisfies the regularity assumption
\begin{equation}
c\in C^{k,\alpha}(\re)\cap C^{\infty}(\re\setminus\{0\}),
\label{th:c-k-a}
\end{equation} 
and the degenerate hyperbolicity assumption 
\begin{equation}
0<c(t)\leq 1
\qquad
\forall t\geq 0.
\label{th:c-dh}
\end{equation}

\item  (Regularity of the solution at initial time) It turns out that
\begin{equation}
(u(0),u'(0))\in\G_{s,r,\beta}(A)\times\G_{s,r,\beta}(A)
\qquad
\forall s>1+\frac{k+\alpha}{2},
\label{th:u0-reg}
\end{equation}
independently of $r>0$ and $\beta\in\re$.

\item (Non-regularity of the solution for all positive times) For every $t>0$ it turns out that
\begin{equation}
(u(t),u'(t))\not\in\G_{-S,R,\beta}(A)\times\G_{-S,R,\beta}(A)
\qquad
\forall S>1+\frac{k+\alpha}{2},
\label{th:ut-schifo}
\end{equation}
independently of $R>0$ and $\beta\in\re$.
\end{enumerate}
	
\end{thm}

\begin{rmk}
\begin{em}

Both Theorem~\ref{thm:main} and Theorem~\ref{thm:dgcs} do not cover the limit case where $\sigma=(2+k+\alpha)^{-1}$. Nevertheless, a careful inspection of the proof of Theorem~\ref{thm:main} reveals that the same conclusions hold true even in the limit case, but provided that $\delta$ is large enough. We skip this more general assumption because it only complicates calculations without introducing new ideas.

On the contrary, the strict inequality in (\ref{hp:sigma-ka-cjs}) seems to be essential in the construction of our counterexamples.

\end{em}
\end{rmk}


\setcounter{equation}{0}
\section{Heuristics}\label{sec:heuristics}

\begin{figure}[htbp]
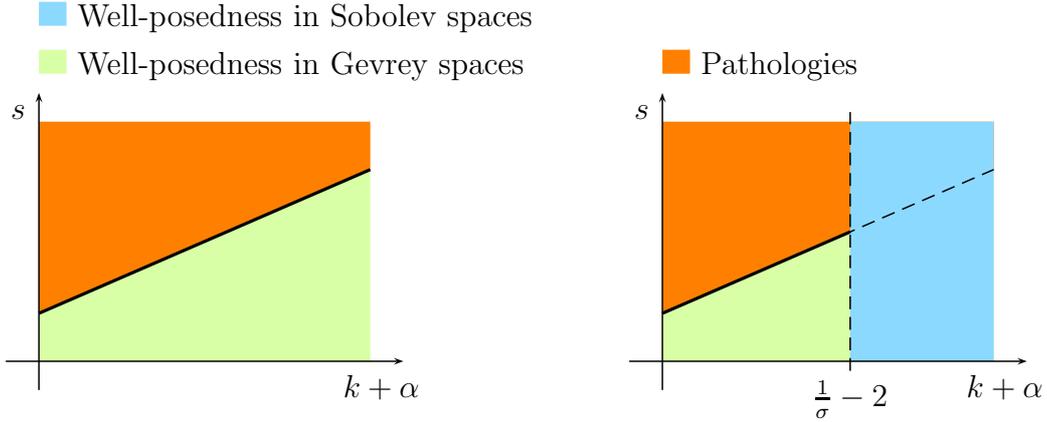

\newrgbcolor{azzurrino}{0.55 0.85 1}
\newrgbcolor{verdino}{0.85 1 0.65}
\psset{xunit=8ex,yunit=7ex}

\noindent
\hfill
\pspicture(-0.5,-0.5)(3.5,4)
\pspolygon*[linecolor=orange](0,0.5)(3,2)(3,2.5)(0,2.5)
\pspolygon*[linecolor=verdino](0,0.5)(3,2)(3,0)(0,0)
\psline[linewidth=0.7\pslinewidth]{->}(-0.3,0)(3.3,0)
\psline[linewidth=0.7\pslinewidth]{->}(0,-0.3)(0,2.8)
\psline[linewidth=1.5\pslinewidth](0,0.5)(3,2)
\uput[-90](3.1,0){$k+\alpha$}
\uput[180](0,2.6){$s$}
\psframe*[linecolor=azzurrino](0,3.5)(0.25,3.75)
\rput[Bl](0.35,3.5){Well-posedness in Sobolev spaces}
\psframe*[linecolor=verdino](0,3)(0.25,3.25)
\rput[Bl](0.35,3){Well-posedness in Gevrey spaces}
\endpspicture
\hfill\hfill\hfill
\pspicture(-0.5,-0.5)(3.5,3.5)
\pspolygon*[linecolor=orange](0,0.5)(3,2)(3,2.5)(0,2.5)
\pspolygon*[linecolor=verdino](0,0.5)(3,2)(3,0)(0,0)
\psframe*[linecolor=azzurrino](1.7,0)(3,2.5)
\psline[linewidth=0.7\pslinewidth]{->}(-0.3,0)(3.3,0)
\psline[linewidth=0.7\pslinewidth]{->}(0,-0.3)(0,2.8)
\psline[linewidth=0.7\pslinewidth, linestyle=dashed](0,0.5)(3,2)
\psline[linewidth=0.7\pslinewidth, linestyle=dashed](1.7,-0.1)(1.7,2.6)
\psline[linewidth=1.5\pslinewidth](0,0.5)(1.7,1.35)
\uput[-90](3.1,0){$k+\alpha$}
\uput[180](0,2.6){$s$}
\rput[t](1.7,-0.2){$\frac{1}{\sigma}-2$}
\psframe*[linecolor=orange](0,3)(0.25,3.25)
\rput[Bl](0.35,3){Pathologies}
\endpspicture
\hfill\mbox{}
\caption{Non-dissipative equation (left) vs dissipative equation (right)}
\label{fig:heuristics}
\end{figure}

The diagrams of Figure~\ref{fig:heuristics} summarize the results of this paper. In the horizontal axis we represent the value $k+\alpha$, corresponding to the time-regularity of $c(t)$. In the vertical axis we represent the space-regularity of initial data, where the value $s$ stands for the Gevrey space of order $s$ (so that higher values of $s$ mean lower regularity).  The oblique line has equation $s=1+\frac{k+\alpha}{2}$. 

For the non-dissipative equation ($\delta=0$) we have the situation described in Theorem~\ref{thmbibl:regularity} and Remark~\ref{rmk:cjs}, namely well-posedness provided that $c(t)$ is of class $C^{k,\alpha}$ and initial data are in Gevrey spaces of order $s<1+\frac{k+\alpha}{2}$, and potential pathologies if $s>1+\frac{k+\alpha}{2}$. The same picture applies if $\delta>0$ and $\sigma=0$.

For the dissipative equation ($\delta>0$) the problem is well-posed \emph{in the Sobolev setting} in the full strip with $k+\alpha>-2+1/\sigma$, as stated in Theorem~\ref{thm:main}.  The region on the left of the vertical line is divided as in the non-dissipative case.  Indeed, Theorem~\ref{thmbibl:regularity} still provides well-posedness \emph{in the Gevrey setting} below the oblique line, while Theorem~\ref{thm:dgcs} shows that pathologies are possible above the oblique line. What happens on the oblique and on the vertical line is less clear, because in these regimes the size of $\delta$ becomes relevant.

Now we present a rough justification of this threshold effect.  As already observed, existence results for problem (\ref{pbm:eqn})--(\ref{pbm:data}) are related to estimates for solutions to the family of ordinary differential equations (\ref{pbm:main-ode}). Let us consider the standard energy function $\mathcal{E}(t):=|\ul'(t)|^{2}+\lambda^{2}|\ul(t)|^{2}$. A classical argument shows that
\begin{equation}
	\mathcal{E}(t)\leq \mathcal{E}(0)
	\exp\left(\lambda t+\lambda\int_{0}^{t}|c(s)|\,ds\right),
	\label{heur:E-dh}
\end{equation}
and this estimate is enough to establish Theorem~\ref{thmbibl:dh}.

If in addition $c(t)$ is of class $C^{k,\alpha}$, and satisfies the degenerate hyperbolicity condition~(\ref{hp:wh}), then (\ref{heur:E-dh}) can be improved to
\begin{equation}
	\mathcal{E}(t)\leq M_{1}\,\mathcal{E}(0)
	\exp\left(M_{2}\lambda^{2/(2+k+\alpha)}t\right)
	\label{heur:E-sh}
\end{equation}
for suitable constants $M_{1}$ and $M_{2}$.  Estimates of this kind are the key point in the proof of both Theorem~\ref{thmbibl:existence} and Theorem~\ref{thmbibl:regularity}.  Moreover,  the pathologies described in Remark~\ref{rmk:cjs} are equivalent to saying that the exponent of $\lambda$ in (\ref{heur:E-sh}) is optimal.

On the other hand, if $\sigma\leq 1/2$ and $c(t)$ is constant, then (\ref{pbm:main-ode}) 
can be explicitly integrated, obtaining that
\begin{equation}
	\mathcal{E}(t)\leq M_{3}\,\mathcal{E}(0)\exp\left(-2\delta\lambda^{2\sigma}t\right)
	\label{heur:E-delta}
\end{equation}
for a suitable constant $M_{3}$.

If $c(t)$ is of class $C^{k,\alpha}$ and satisfies the degenerate hyperbolicity condition~(\ref{hp:wh}), then we expect a superposition of the effects of the coefficient, represented by (\ref{heur:E-sh}), and the effects of the damping, represented by (\ref{heur:E-delta}). We end up with something like
\begin{equation}
	\mathcal{E}(t)\leq M_{1}M_{3}\,\mathcal{E}(0)\exp\left((M_{2}\lambda^{2/(2+k+\alpha)}-2\delta\lambda^{2\sigma})\cdot t\right).
	\label{heur:E-conflict}
\end{equation}

Therefore, it is reasonable to expect well-posedness in Sobolev spaces when the argument of the exponential is bounded from above independently of $\lambda$, which is true for sure when condition (\ref{hp:sigma-ka}) is satisfied.  On the contrary, when (\ref{hp:sigma-ka-cjs}) is satisfied, the right-hand side of (\ref{heur:E-conflict}) diverges as $\lambda\to+\infty$, opening the door to the pathologies. In the border-line case, namely when the two exponents are equal, the size of $\delta$ comes into play. 


\setcounter{equation}{0}
\section{Proofs of well-posedness and regularity results}\label{sec:proofs}

In this section we prove Theorem~\ref{thm:main}. The proof has three main steps.
\begin{itemize}

\item  In Lemma~\ref{lemma:gamma} we show that $c(t)$ can be approximated by a family $\gl(t)$ of nonnegative functions of class $C^{1}$ satisfying suitable estimates. Glaeser type inequalities play a crucial role in this step.

\item In Proposition~\ref{prop:ODE} we use the functions $\gl(t)$ as coefficients of approximated hyperbolic energies, and with the help of these energies we estimate the growth of solutions to the family of ordinary differential equations (\ref{pbm:main-ode}).

\item Finally, we conclude by means of the spectral theory and the previous estimates.

\end{itemize}

In the sequel we set for simplicity
\begin{equation}
\theta:=\frac{2}{2+k+\alpha},
\label{defn:theta}
\end{equation}
and we observe that definition (\ref{defn:theta}) implies that
\begin{equation}
\frac{2(1-\theta)}{k+\alpha}=\theta.
\label{defn:theta-bis}
\end{equation}

\begin{lemma}[Approximation of the coefficient]\label{lemma:gamma}

Let $T>0$, $\mu>0$ and  $\alpha\in(0,1]$ be real numbers, let $k\in\n$ be a nonnegative integer, and let $c:[0,T]\to[0,\mu]$ be any function of class $C^{k,\alpha}$. Let $\theta$ be defined by (\ref{defn:theta}).

Then for every $\lambda>0$ there exists a function $\gl:[0,T]\to[0,\mu]$ of classe $C^{1}$ such that
\begin{equation}
|c(t)-\gl(t)|\leq\lambda^{-2(1-\theta)}
\qquad
\forall t\in[0,T],
\label{th:c-g}
\end{equation}
and
\begin{equation}
|c(t)-\gl(t)|^{2}\leq c(t)\lambda^{-2(1-\theta)}
\qquad
\forall t\in[0,T].
\label{th:c-g-2}
\end{equation}

In addition, the derivative of $\gl$ satisfies the following estimates depending on $k$.\begin{itemize}
  \item If $k=0$, then it turns out that
  \begin{equation}
  |\gl'(t)|\leq (25H)^{1/\alpha}\cdot c(t)\lambda^{\theta}
\qquad
\forall t\in[0,T],
  \label{th:gl'-0}
  \end{equation}
  where $H$ denotes the $\alpha$-H\"older constant of $c(t)$ in $[0,T]$.
  
  \item If $k=1$, then it turns out that
  \begin{equation}
  |\gl'(t)|\leq 4K\cdot c(t)\lambda^{\theta}
\qquad
\forall t\in[0,T],
  \label{th:gl'-1}
  \end{equation}
  where $K$ denotes the constant for which (\ref{th:glaeser-1}) holds true.
  
  \item If $k\geq 2$, then it turns out that
  \begin{equation}
  |\gl'(t)|\leq 4\varphi(t)\cdot c(t)\lambda^{\theta}
\qquad
\forall t\in[0,T],
  \label{th:gl'-2}
  \end{equation}
  where $\varphi(t)$ denotes the function for which (\ref{th:glaeser-2}) holds true.

\end{itemize}

\end{lemma}

\paragraph{\textmd{\textit{Proof}}}

For every $\ep>0$, let us consider the function $\psi_{\ep}:\re\to\re$ defined as
\begin{equation}
\psi_{\ep}(x):=\left\{
\begin{array}{l@{\quad}l}
0  &  \mbox{if } x\leq 0,  \\[1ex]
\displaystyle x-\frac{2}{\pi}\cdot\ep\arctan\left(\frac{\pi}{2}\cdot\frac{x}{\ep}\right)  &  \mbox{if } x\geq 0.     
\end{array}
\right.
\nonumber
\end{equation}

It turns out that $\psi_{\ep}$ is a function of class $C^{1}$ that approximates the piecewise affine function $\max\{x,0\}$. In particular, in the sequel we need that
\begin{equation}
|\psi_{\ep}(x)-x|\leq\ep
\qquad
\forall\ep>0,\quad\forall x\geq 0,
\label{psiep-x}
\end{equation}
and that the derivative satisfies
\begin{equation}
\psi_{\ep}'(x)=0
\qquad
\forall\ep>0,\quad\forall x\leq 0,
\label{psiep'=0}
\end{equation}
\begin{equation}
|\psi_{\ep}'(x)|\leq 1
\qquad
\forall\ep>0,\quad\forall x\geq 0.
\label{psiep'<1}
\end{equation}

\subparagraph{\textmd{\textit{Case $k\geq 1$}}}

In this case we set
\begin{equation}
\cep(t):=\psi_{\ep}(c(t)-\ep)
\qquad
\forall t\in[0,T],
\nonumber
\end{equation}
and we define $\gl(t)$ as $\cep(t)$ in the case where
\begin{equation}
4\ep:=\lambda^{-2(1-\theta)}.
\label{defn:ep-k>=1}
\end{equation}

In order to prove (\ref{th:c-g}) and (\ref{th:c-g-2}), we distinguish two cases. 
\begin{itemize}

\item  If $c(t)\leq\ep$, then $\cep(t)=0$, and hence
$$|\gl(t)-c(t)|=c(t)\leq\ep\leq\lambda^{-2(1-\theta)},$$
which proves (\ref{th:c-g}). Similarly, it turns out that
$$|\gl(t)-c(t)|^{2}=c(t)\cdot c(t)\leq c(t)\cdot\ep\leq c(t)\cdot\lambda^{-2(1-\theta)},$$
which proves (\ref{th:c-g-2}).

\item  If $c(t)\geq\ep$, then from (\ref{psiep-x}) with $x=c(t)-\ep$ we deduce that
$$|\gl(t)-c(t)|=|\psi_{\ep}(c(t)-\ep)-(c(t)-\ep)-\ep|\leq 2\ep\leq\lambda^{-2(1-\theta)},$$
which proves (\ref{th:c-g}), and similarly
$$|\gl(t)-c(t)|^{2}\leq 2\ep\cdot 2\ep\leq c(t)\cdot  4\ep=c(t)\cdot\lambda^{-2(1-\theta)},$$
which proves (\ref{th:c-g-2}).

\end{itemize}

Let us consider now the derivative $\gl'(t)$, and let us distinguish again two cases.
\begin{itemize}

\item  If $c(t)\leq\ep$, then from (\ref{psiep'=0}) we deduce that $\gl'(t)=0$, and hence both (\ref{th:gl'-1}) and (\ref{th:gl'-2}) are trivial.

\item  If $c(t)\geq\ep$, then from (\ref{psiep'<1}) we deduce that
\begin{equation}
|\gl'(t)|=|\psi_{\ep}'(c(t)-\ep)|\cdot|c'(t)|\leq|c'(t)|.
\label{gl'-c}
\end{equation}

Now we apply the Glaeser type inequalities of Theorem~\ref{thmbibl:glaeser}. If $k=1$, from (\ref{defn:ep-k>=1}) and (\ref{defn:theta-bis}) we obtain that
$$|c'(t)|\leq Kc(t)\cdot[c(t)]^{-1/(1+\alpha)}\leq Kc(t)\ep^{-1/(1+\alpha)}\leq 4Kc(t)\lambda^{\theta}.$$

Plugging this inequality into (\ref{gl'-c}) we obtain (\ref{th:gl'-1}). Similarly, for $k\geq 2$ we obtain that
$$|c'(t)|\leq\varphi(t)c(t)\cdot[c(t)]^{-1/(k+\alpha)}\leq\varphi(t)c(t)\ep^{-1/(k+\alpha)}\leq 4\varphi(t)c(t)\lambda^{\theta}.$$

Plugging this inequality into (\ref{gl'-c}) we obtain (\ref{th:gl'-2}).

\end{itemize}  

\subparagraph{\textmd{\textit{Case $k=0$}}}

In this case we first extend $c(t)$ to the whole half-line $t\geq 0$ by setting $c(t)=c(T)$ for every $t\geq T$, and then we consider the regularized function
$$\ceph(t):=\frac{1}{\ep}\int_{t}^{t+\ep}c(s)\,ds
\qquad
\forall t\geq 0.$$

Due to (\ref{hp:wh}), the function $\ceph(t)$ takes its values in $[0,\mu]$. Moreover, it is of class $C^{1}$ and satisfies
\begin{equation}
|\ceph(t)-c(t)|\leq H\ep^{\alpha}
\qquad
\forall t\in[0,T],
\label{est:ceph}
\end{equation}
and
\begin{equation}
|\ceph\mbox{}\!'(t)|\leq\frac{H}{\ep^{1-\alpha}}
\qquad
\forall t\in[0,T].
\label{est:ceph'}
\end{equation}

Now we set
\begin{equation}
\cep(t):=\psi_{2H\ep^{\alpha}}(\ceph(t)-2H\ep^{\alpha})
\qquad
\forall t\in[0,T],
\nonumber
\end{equation}
and we define $\gl(t)$ as $\cep(t)$ in the case where
\begin{equation}
25H\ep^{\alpha}:=\lambda^{-2(1-\theta)}.
\label{defn:ep-k=0}
\end{equation}

In order to prove (\ref{th:c-g}) and (\ref{th:c-g-2}), we distinguish two cases. 
\begin{itemize}

\item  If $\ceph(t)\leq 2H\ep^{\alpha}$, then $\cep(t)=0$, and in addition $c(t)\leq 3H\ep^{\alpha}$ because of (\ref{est:ceph}) and triangle inequality. It follows that
$$|\gl(t)-c(t)|=c(t)\leq 3H\ep^{\alpha}\leq\lambda^{-2(1-\theta)},$$
which proves (\ref{th:c-g}). Similarly, we obtain that
$$|\gl(t)-c(t)|^{2}=c(t)\cdot c(t)\leq c(t)\cdot 3H\ep^{\alpha}\leq c(t)\cdot\lambda^{-2(1-\theta)},$$
which proves (\ref{th:c-g-2}).

\item  If $\ceph(t)\geq 2H\ep^{\alpha}$, then $c(t)\geq H\ep^{\alpha}$ because of (\ref{est:ceph}) and triangle inequality. Now we observe that
$$|\gl(t)-c(t)|=\left|\psi_{2H\ep^{\alpha}}(\ceph(t)-2H\ep^{\alpha})-(\ceph(t)-2H\ep^{\alpha})+(\ceph(t)-c(t))-2H\ep^{\alpha}\strut\right|.$$

The first two terms can be estimated by means of inequality (\ref{psiep-x}) with $x=\ceph(t)-2H\ep^{\alpha}$. The third term can be estimates as in (\ref{est:ceph}). Thus from triangle inequality we deduce that
$$|\gl(t)-c(t)|\leq 5H\ep^{\alpha}\leq\lambda^{-2(1-\theta)},$$
and
$$|\gl(t)-c(t)|^{2}\leq 5H\ep^{\alpha}\cdot 5H\ep^{\alpha}\leq c(t)\cdot 25H\ep^{\alpha}=c(t)\cdot\lambda^{-2(1-\theta)},$$
which prove (\ref{th:c-g}) and (\ref{th:c-g-2}).

\end{itemize}

As for derivatives, again we distinguish two cases. 
\begin{itemize}

\item  If $c(t)\leq H\ep^{\alpha}$, then from (\ref{est:ceph}) we deduce that $\ceph(t)\leq 2H\ep^{\alpha}$. Thus from (\ref{psiep'=0}) we conclude that $\gl'(t)=0$, and hence (\ref{th:gl'-0}) is trivial. 

\item  If $c(t)\geq H\ep^{\alpha}$, then from (\ref{psiep'=0}) and (\ref{psiep'<1}) we deduce that
\begin{equation}
|\gl'(t)|=\left|\psi_{2H\ep^{\alpha}}'(\ceph(t)-2H\ep^{\alpha})\right|\cdot|\ceph\mbox{}\!'(t)|\leq|\ceph\mbox{}\!'(t)|,
\nonumber
\end{equation}
and therefore from (\ref{est:ceph'}), (\ref{defn:ep-k=0}) and (\ref{defn:theta-bis}) we conclude that
$$|\gl'(t)|\leq|\ceph\mbox{}\!'(t)|\leq\frac{H\ep^{\alpha}}{\ep}\leq c(t)\cdot\frac{1}{\ep}=c(t)\cdot(25H)^{1/\alpha}\lambda^{\theta},$$
which proves (\ref{th:gl'-0}).

\end{itemize}

This completes the proof.\qed

\begin{prop}[Estimates on components]\label{prop:ODE}
	
Let us consider problem (\ref{pbm:main-ode})--(\ref{eqn:ODE-data}) under the following assumptions:
\begin{itemize}
		
\item  there exists $k\in\n$ and $\alpha\in(0,1]$ such that $c\in\cka([0,T])$,
	
\item  $c(t)$ satisfies the degenerate hyperbolicity assumption (\ref{hp:wh}) and condition (\ref{hp:c0=0}),

\item $\delta$ and $\lambda$ are positive real numbers, and $\sigma$ is a real number satisfying (\ref{hp:sigma-ka}).
	
\end{itemize}

Then there exist positive real numbers $r$ and $\nu$, both independent of $\lambda$, such that the following estimates hold true.
\begin{enumerate}
	\renewcommand{\labelenumi}{(\arabic{enumi})}
		
\item \emph{(Case $k\in\{0,1\}$)} For every $t\in[0,T]$, and every $\lambda\geq\nu$, it turns out that
\begin{equation}
	|\ul'(t)|^{2}+\delta^{2}\lambda^{4\sigma}|\ul(t)|^{2\sigma}\leq 3\left(u_{1,\lambda}^{2}+\delta^{2}\lambda^{4\sigma}u_{0,\lambda}^{2}\right)\exp\left(-4r\lambda^{2\sigma}C(t)\right),
	\label{th:est-u-01}
\end{equation}
where $C(t)$ is defined by (\ref{defn:C}).
			
\item \emph{(Case $k\geq 2$)} Let $\varphi(t)$ be the function which appears in (\ref{th:glaeser-2}), let $\theta$ be defined by (\ref{defn:theta}), and let
$$\Phi(t):=\int_{0}^{t}\varphi(s)\,ds
\qquad
\forall t\in [0,T].$$

Then for every $t\in[0,T]$, and every $\lambda\geq\nu$, it turns out that
\begin{eqnarray}
|\ul'(t)|^{2}+\delta^{2}\lambda^{4\sigma}|\ul(t)|^{2\sigma} & \leq & 3\left(u_{1,\lambda}^{2}+\delta^{2}\lambda^{4\sigma}u_{0,\lambda}^{2}\right)\cdot\mbox{}  
\nonumber \\[1ex]
 & & \mbox{}\cdot\exp\left(-4r\lambda^{2\sigma}C(t)+4\lambda^{\theta}\Phi(t)\right).
 \label{th:est-u-2}
\end{eqnarray}
			
\end{enumerate}
\end{prop}

\paragraph{\textmd{\textit{Proof}}}

For every function $\gamma:[0,T]\to[0,+\infty)$ of class $C^{1}$, we introduce the approximated hyperbolic energy
\begin{equation}
\eg(t):=|\ul'(t)|^{2}+\delta^{2}\lambda^{4\sigma}|\ul(t)|^{2}+\delta\lambda^{2\sigma}\ul(t)\ul'(t)+\gamma(t)\lambda^{2}|\ul(t)|^{2}.
\nonumber
\end{equation}

Since
$$\delta\lambda^{2\sigma}|\ul(t)\ul'(t)|\leq\frac{\delta^{2}}{2}\lambda^{4\sigma}|\ul(t)|^{2}+\frac{1}{2}|\ul'(t)|^{2},$$
it follows that
\begin{equation}
\eg(t)\geq\frac{1}{2}|\ul'(t)|^{2}+\frac{\delta^{2}}{2}\lambda^{4\sigma}|\ul(t)|^{2}+\gamma(t)\lambda^{2}|\ul(t)|^{2}
\label{est:eg>}
\end{equation}
and
\begin{equation}
\eg(t)\leq\frac{3}{2}|\ul'(t)|^{2}+\frac{3\delta^{2}}{2}\lambda^{4\sigma}|\ul(t)|^{2}+\gamma(t)\lambda^{2}|\ul(t)|^{2}
\label{est:eg<}
\end{equation}
for every admissible value of the parameters.

Given any real number $r$, an elementary but lengthy calculation shows that
\begin{equation}
\eg'(t)=-4r c(t)\lambda^{2\sigma}\eg(t)-Q_{1,\gamma}(t)+Q_{2,\gamma}(t)
\qquad
\forall t\in[0,T],
\label{eqn:eg'}
\end{equation}
where
$$Q_{1,\gamma}(t):=X_{\gamma}(t)|\ul'(t)|^{2}+Y_{\gamma}(t)|\ul(t)|^{2}+Z_{\gamma}(t)\ul(t)\ul'(t)$$
is a quadratic form in the variables $\ul(t)$ and $\ul'(t)$ with coefficients
$$X_{\gamma}(t):=\lambda^{2\sigma}(3\delta-4rc(t)),
\hspace{3em}
Y_{\gamma}(t):= c(t)\lambda^{2+2\sigma}\left(\frac{\delta}{2}-4r\delta^{2}\lambda^{4\sigma-2}-4r\gamma(t)\right),$$
$$Z_{\gamma}(t):=2(c(t)-\gamma(t))\lambda^{2}-4r\delta c(t)\lambda^{4\sigma},$$
and
$$Q_{2,\gamma}(t):=-\frac{\delta}{2}c(t)\lambda^{2+2\sigma}|\ul(t)|^{2}+\gamma'(t)\lambda^{2}|\ul(t)|^{2}.$$

In the sequel we fix $r$ such that
\begin{equation}
2r\mu\leq\delta,
\hspace{3em}
16r(\delta^{2}+\mu)\leq\delta,
\hspace{3em}
64r^{2}\mu\leq 1,
\label{defn:r}
\end{equation}
and we provide estimates for $Q_{1,\gamma}(t)$ and $Q_{2,\gamma}(t)$.

\subparagraph{\textmd{\textit{Estimate on $Q_{1,\gamma}(t)$}}}

Let $\nu$ be a positive real number such that
\begin{equation}
\nu\geq 1,
\hspace{4em}
\delta\nu^{2\sigma-\theta}\geq 4.
\label{defn:nu}
\end{equation}

For every $\lambda>0$, let $\gl(t)$ be the approximation of $c(t)$ provided by Lemma~\ref{lemma:gamma}. We claim that
\begin{equation}
Q_{1,\gl}(t)\geq 0
\qquad
\forall \lambda\geq\nu,\quad\forall t\in[0,T].
\label{est:Q1>0}
\end{equation}

From the theory of quadratic forms, we know that (\ref{est:Q1>0}) holds true if the three inequalities
\begin{equation}
X_{\gl}(t)\geq 0,
\hspace{3em}
Y_{\gl}(t)\geq 0,
\hspace{3em}
4X_{\gl}(t)Y_{\gl}(t)\geq Z_{\gl}(t)^{2}
\label{est:XY>ZZ}
\end{equation}
are satisfied for every $\lambda\geq\nu$ and every $t\in[0,T]$.

From (\ref{hp:wh}) and the first inequality in (\ref{defn:r}) we obtain that $X_{\gl}(t)\geq\delta\lambda^{2\sigma}$. Moreover, since $\sigma\leq 1/2$, $\lambda\geq 1$, and $\gl(t)\leq\mu$, from the second inequality in (\ref{defn:r}) we obtain that $Y_{\gl}\geq\delta c(t)\lambda^{2+2\sigma}/4$. This proves the first two inequalities in (\ref{est:XY>ZZ}), and also provides the following estimate
\begin{equation}
4X_{\gl}(t)Y_{\gl}(t)\geq\delta^{2}c(t)\lambda^{2+4\sigma}
\label{est:XY}
\end{equation}
for the left-hand side of  the third one. As for the right-hand side, we first observe that
\begin{equation}
Z_{\gl}(t)^{2}\leq 8(c(t)-\gl(t))^{2}\lambda^{4}+32r^{2}\delta^{2}c(t)^{2}\lambda^{8\sigma}.
\label{est:ZZ-1}
\end{equation}

The second term can be estimated as
\begin{equation}
32r^{2}\delta^{2}c(t)^{2}\lambda^{8\sigma}\leq 32r^{2}\mu\lambda^{4\sigma-2}\cdot\delta^{2}c(t)\lambda^{2+4\sigma}\leq\frac{\delta^{2}}{2}c(t)\lambda^{2+4\sigma},
\nonumber
\end{equation}
where again we used that $\sigma\leq 1/2$, $\lambda\geq 1$, and the last inequality in (\ref{defn:r}). As for the first term, now we exploit the special choice of $\gl(t)$ provided by Lemma~\ref{lemma:gamma}. From (\ref{th:c-g-2}) and (\ref{defn:nu}) we obtain that
\begin{equation}
8(c(t)-\gl(t))^{2}\lambda^{4}\leq 8c(t)\lambda^{-2(1-\theta)}\lambda^{4}\leq\frac{\delta^{2}}{2}c(t)\lambda^{2+4\sigma}.
\nonumber
\end{equation}

Plugging the last two estimates into (\ref{est:ZZ-1}) we conclude that
\begin{equation}
Z_{\gl}(t)^{2}\leq\delta^{2}c(t)\lambda^{2+4\sigma}.
\label{est:ZZ}
\end{equation}

Finally, from (\ref{est:XY}) and (\ref{est:ZZ}) we obtain the third inequality in (\ref{est:XY>ZZ}), and this completes the proof of (\ref{est:Q1>0}).

\subparagraph{\textmd{\textit{Estimate on $Q_{2,\gamma}(t)$ and conclusion if $k=0$}}}

Let $H$ denote the $\alpha$-H\"older constant of $c(t)$ in $[0,T]$, and let us assume that $\nu$ satisfies (\ref{defn:nu}) and the further condition
\begin{equation}
\nu^{2\sigma-\theta}\geq\frac{2}{\delta}(25H)^{1/\alpha}.
\label{hp:nu-k0}
\end{equation}

As before, let $\gl(t)$ denote the approximation of $c(t)$ provided by Lemma~\ref{lemma:gamma}. From (\ref{th:gl'-0}) and (\ref{hp:nu-k0}) it follows that
\begin{equation}
Q_{2,\gl}(t)\leq 0
\qquad
\forall \lambda\geq\nu,\quad\forall t\in[0,T].
\label{est:Q2-k0}
\end{equation}

Plugging (\ref{est:Q1>0}) and (\ref{est:Q2-k0}) into (\ref{eqn:eg'}), we find that
$$\egl'(t)\leq -4r c(t)\lambda^{2\sigma}\egl(t)
\qquad
\forall \lambda\geq\nu,\quad\forall t\in[0,T].$$

Integrating this differential inequality we obtain that
\begin{equation}
\egl(t)\leq\egl(0)\exp\left(-4r C(t)\lambda^{2\sigma}\right)
\qquad
\forall \lambda\geq\nu,\quad\forall t\in[0,T].
\label{est:eg}
\end{equation}

Finally, we observe that $\gl(0)=0$ because of (\ref{th:c-g-2}) and our assumption that $c(0)=0$. At this point, estimate (\ref{th:est-u-01}) with $k=0$ follows from (\ref{est:eg>}), (\ref{est:eg<}), and (\ref{est:eg}). 

\subparagraph{\textmd{\textit{Estimate on $Q_{2,\gamma}(t)$ and conclusion if $k=1$}}}

Let $K$ denote the constant such that (\ref{th:glaeser-1}) holds true, and let us assume that $\nu$ satisfies (\ref{defn:nu}) and the further condition
\begin{equation}
\nu^{2\sigma-\theta}\geq\frac{8K}{\delta}.
\label{hp:nu-k1}
\end{equation}

As before, let $\gl(t)$ denote the approximation of $c(t)$ provided by Lemma~\ref{lemma:gamma}. From (\ref{th:gl'-1}) and (\ref{hp:nu-k1}) it follows that also in this case (\ref{est:Q2-k0}) holds true.

At this point, the conclusion follows exactly as in the case $k=0$.

\subparagraph{\textmd{\textit{Estimate on $Q_{2,\gamma}(t)$ and conclusion if $k\geq 2$}}}

Let us assume that $\nu$ satisfies (\ref{defn:nu}) and the further condition
\begin{equation}
\nu^{2\sigma-\theta}\geq\frac{\sqrt{2}}{\delta}.
\label{hp:nu-k2}
\end{equation}

As always, let $\gl(t)$ denote the approximation of $c(t)$ provided by Lemma~\ref{lemma:gamma}. From (\ref{th:gl'-2}) it follows that
$$|\gl'(t)| \leq 4c(t)\varphi(t)\lambda^{\theta} = 4\gl(t)\varphi(t)\lambda^{\theta}+4(c(t)-\gl(t))\varphi(t)\lambda^{\theta}.$$

On the other hand, from (\ref{th:c-g}) and (\ref{hp:nu-k2}) it follows that
$$|c(t)-\gl(t)|\leq\lambda^{-2(1-\theta)}\leq\frac{\delta^{2}}{2}\lambda^{4\sigma-2},$$
and hence from (\ref{est:eg>}) we deduce that
\begin{eqnarray*}
Q_{2,\gl}(t)  &  \leq & \lambda^{2}|\gl'(t)|\cdot|\ul(t)|^{2}  \\[1ex]
  & \leq &  4\lambda^{\theta}\varphi(t)\left[\gl(t)\lambda^{2}|\ul(t)|^{2}+\frac{\delta^{2}}{2}\lambda^{4\sigma}|\ul(t)|^{2}\right]  \\[1ex]
  & \leq &  4\lambda^{\theta}\varphi(t)\egl(t).
\end{eqnarray*}

Plugging this estimate and (\ref{est:Q1>0}) into (\ref{eqn:eg'}), we find that
$$\egl'(t)\leq \left[-4r c(t)\lambda^{2\sigma}+4\lambda^{\theta}\varphi(t)\right]\egl(t).$$

Integrating this differential inequality we obtain that
$$\egl(t)\leq\egl(0)\exp\left(-4r C(t)\lambda^{2\sigma}+4\lambda^{\theta}\Phi(t)\right)$$
for every $\lambda\geq\nu$ and every $t\in[0,T]$. Recalling that $\gl(0)=0$, the conclusion (\ref{th:est-u-2}) follows again from (\ref{est:eg>}) and (\ref{est:eg<}), as in the previous cases.\qed

\begin{rmk}
\begin{em}

We observe that $r$ depends only on $\delta$ and $\mu$, as specified by (\ref{defn:r}), and that $\nu$ depends on $\delta$, $\sigma$, $k+\alpha$, and
\begin{itemize}
\item on the $\alpha$-H\"older constant of $c(t)$ when $k=0$,
\item  on the constant $K$ for which (\ref{th:glaeser-1}) is true when $k=1$.
\end{itemize}

The conditions on $\nu$ are stated in (\ref{defn:nu}), and in (\ref{hp:nu-k0}), (\ref{hp:nu-k1}) or (\ref{hp:nu-k2}), depending on the value of $k$.

Stating precisely ``what depends on what'' could be useful when considering families of equations of the form (\ref{pbm:eqn}) with different choices of $\delta$, $\sigma$, $c(t)$. This is often a key step in the fixed point arguments exploited when dealing with nonlinear problems.

\end{em}
\end{rmk}


\subsubsection*{End of the proof of Theorem~\ref{thm:main}}

Let $r$ and $\nu$ be as in Proposition~\ref{prop:ODE}. Let us write $H$ as an orthogonal direct sum
\begin{equation}
	H:=H_{\nu,-}\oplus H_{\nu,+},
\nonumber
\end{equation}
where $H_{\nu,-}$ is the closure of the subspace generated by all eigenvectors of $A$ relative to eigenvalues $\lambda_{n}<\nu$, and $H_{\nu,+}$ is the closure of the subspace generated by all eigenvectors of $A$ relative to eigenvalues $\lambda_{n}\geq\nu$. Let $u_{\nu,-}(t)$ and $u_{\nu,+}(t)$ denote the corresponding components of $u(t)$.

The low frequency component $u_{\nu,-}(t)$ is continuous in any reasonable space because the operator $A$ is bounded in $H_{\nu,-}$. For further details we refer to~\cite[Remark~3.3]{gg:dgcs-strong}.

In order to estimate the high frequency component $u_{\nu,+}(t)$, we apply Proposition~\ref{prop:ODE} to all components $u_{n}(t)$ of $u(t)$ corresponding to eigenvalues $\li\geq\nu$. To this end, we distinguish two cases.

\subparagraph{\textmd{\textit{Case $k\in\{0,1\}$}}}

For these values of $k$ we know that estimate (\ref{th:est-u-01}) holds true for every $t\in[0,T]$. Summing over all eigenvalues $\li\geq\nu$, we obtain that $u_{\nu,+}(t)$ is bounded in $\G_{(2\sigma)^{-1},rC(t),\sigma}(A)$, and hence in particular in $D(A^{\sigma})$, while $u_{\nu,+}'(t)$ is bounded in $\G_{(2\sigma)^{-1},rC(t),0}(A)$, and hence in particular in $H$.  

The same estimate guarantees the uniform convergence in $[0,T]$ of the series defining $u_{\nu,+}(t)$ and $u_{\nu,+}'(t)$ in the same spaces.  Since all summands are continuous, and the convergence is uniform, the sum is continuous as well. 

This proves (\ref{th:main-sob-open}) through (\ref{th:main-sob-closed}) for these values of $k$.

\subparagraph{\textmd{\textit{Case $k\geq 2$}}}

For these values of $k$ we have to rely on estimate (\ref{th:est-u-2}), which is worse than (\ref{th:est-u-01}) because $C(t)$ vanishes of order at least one in $t=0$, while $\Phi(t)$ might vanish with a lower order.

Nevertheless, if we fix any $\tau\in(0,T)$, from assumption (\ref{hp:C>0}) we deduce that $C(t)$ is bounded from below in $[\tau,T]$ by a positive constant, while $\Phi(t)$ is bounded from above in the same interval. Recalling (\ref{hp:sigma-ka}), we obtain that
$$-4r\lambda^{2\sigma}C(t)+4\lambda^{\theta}\Phi(t)\leq-2r\lambda^{2\sigma}C(t)
\qquad
\forall t\in[\tau,T]$$
provided that $\lambda$ is large enough.

At this point, the same argument of the previous case proves that 
$$u\in C^{0}\left([\tau,T],\G_{(2\sigma)^{-1},rC(t),\sigma}(A)\right)\cap C^{1}\left([\tau,T],\G_{(2\sigma)^{-1},rC(t),0}(A)\right).$$

Since $\tau$ is arbitrary, this proves (\ref{th:main-sob-open}) and (\ref{th:main-gevrey}) for these values of $k$.\qed


\begin{rmk}\label{rmk:estimates}
\begin{em}

The proof of Theorem~\ref{thm:main} provides also decay estimates in Sobolev spaces for high frequency components of solutions. Indeed, in the case $k\in\{0,1\}$ from (\ref{th:est-u-01}) we obtain that
\begin{equation}
|u_{\nu,+}'(t)|^{2}+\delta^{2}|A^{\sigma}u_{\nu,+}(t)|^{2}\leq 3\left(|u_{\nu,+}'(0)|^{2}+\delta^{2}|A^{\sigma}u_{\nu,+}(0)|^{2}\right)\exp\left(-4r\nu^{2\sigma}C(t)\right).
\nonumber
\end{equation}

Analogous estimates can be deduced from (\ref{th:est-u-2}) in the case $k\geq 2$.

On the contrary, low frequency components could even grow linearly with time (it is enough to think to the case $\lambda=0$).

\end{em}
\end{rmk}

\begin{rmk}\label{rmk:ka}
\begin{em}

As announced at the beginning of the section, the proof of Theorem~\ref{thm:main} followed the  path
$$
\fbox{$u\in C^{k,\alpha}$}
\Rightarrow
\fbox{Glaeser}
\Rightarrow
\fbox{Lemma~\ref{lemma:gamma}}
\Rightarrow
\fbox{Prop.~\ref{prop:ODE}}
\Rightarrow
\fbox{Theorem~\ref{thm:main}}\,,
$$
where each step depends only on the previous one. As a consequence, Theorem~\ref{thm:main} holds true whenever $c(t)$ can be approximated as in Lemma~\ref{lemma:gamma}. More precisely, statements~(1) and~(2) of Theorem~\ref{thm:main} are true for every $\sigma$ satisfying (\ref{hp:sigma-ka}) provided that $c(t)$ can be approximated by a family $c_{\lambda}(t)$ of nonnegative functions of class $C^{1}$ satisfying (\ref{th:c-g}), (\ref{th:c-g-2}), and (\ref{th:gl'-2}) for some $\varphi\in L^{1}((0,T))$. If in addition $\varphi\in L^{\infty}((0,T))$, then also statement~(3) of Theorem~\ref{thm:main} is true.

Just to give an example, let us consider the coefficient $c(t):=t^{30}$. Since this coefficient is of class $C^{\infty}$, from Theorem~\ref{thm:main} we deduce that a solution in $D(A^{\sigma})\times H$ exists for every $\sigma\in(0,1/2]$ (and also for $\sigma>1/2$, as already proved in~\cite{gg:dgcs-strong}). If we want continuity in the same space up to $t=0$, and we limit ourselves to  Theorem~\ref{thm:main} as stated, we can consider $t^{30}$ as a coefficient of class $C^{1,1}$, so that from (\ref{th:main-sob-closed}) we deduce the required continuity in the range $\sigma\in(1/4,1/2]$. 

On the other hand, the function $t^{30}$ satisfies the Glaeser type inequality (\ref{th:glaeser-2}) with $\varphi(t)$ bounded even if $k=29$ and $\alpha=1$. Therefore, we can obtain the conclusions of Lemma~\ref{lemma:gamma} also for this values of $k$ and $\alpha$, and hence deduce that (\ref{th:main-sob-closed}) is valid even in the larger range $\sigma\in(1/32,1/2]$. 

What happens when $\sigma\leq 1/32$ remains unclear.

\end{em}
\end{rmk}


\setcounter{equation}{0}
\section{Counterexamples}\label{sec:counterexamples}

In this section we prove Theorem~\ref{thm:dgcs}. The strategy of the proof dates back to \cite{dgcs}, but the functions and the sequences involved are different from case to case.

The starting point is finding a family of $\lambda$-dependent coefficients for which the ordinary differential equation (\ref{pbm:main-ode}) admits solutions whose energy grows exponentially with time. Then we glue together these $\lambda$-dependent coefficients in order to produce a unique $\lambda$-independent coefficient $c(t)$ that acts on infinitely many time-scales, and realizes a similar growth for countably many components. To this end, we introduce a suitable decreasing sequence $t_{n}\to 0^{+}$, and in the interval $[t_{n},t_{n-1}]$ we design $c(t)$ so that $u_{n}(t_{n})$ is small and $u_{n}(t_{n-1})$ is huge. Then we check that the piecewise defined coefficient $c(t)$ has the required time-regularity, and that $u_{n}(t)$ remains small for $t\in[0,t_{n}]$ and huge for $t\geq t_{n-1}$. This completes the proof.

\subparagraph{\textmd{\textit{Basic ingredients}}}

Let us consider the functions
\begin{equation}
b(m,\ep,\lambda,t):=(2m\lambda\ep-\delta\lambda^{2\sigma})t-\ep\sin(2m\lambda t),
\label{defn:belt}
\end{equation}
\begin{equation}
	w(m,\ep,\lambda,t):=\sin(m\lambda t)\exp(b(m,\ep,\lambda,t)),
	\label{defn:welt}
\end{equation}
\begin{equation}
	\gamma(m,\ep,\lambda,t):=m^{2}+\frac{\delta^{2}}{\lambda^{2-4\sigma}}-16m^{2}\ep^{2}\sin^{4}(m\lambda t)-8m^{2}\ep\sin(2m\lambda t).
	\label{defn:gelt}
\end{equation}

These functions depend on four variables, and generalize the corresponding functions of three variables introduced in~\cite{gg:dgcs-strong}. The key point is that for every admissible value of the parameters it turns out that
$$w''(m,\ep,\lambda,t)+2\delta\lambda^{2\sigma} w'(m,\ep,\lambda,t)+\lambda^{2}\gamma(m,\ep,\lambda,t)w(m,\ep,\lambda,t)=0,$$
where ``primes'' denote differentiation with respect to $t$. If we set parameters in such a way that
$$m\lambda\ep\sim\lambda^{\theta},$$
with $\theta$ defined by (\ref{defn:theta}), then the powers of $\lambda$ in the coefficient of the linear term in (\ref{defn:belt}) are exactly the same powers that appear in the argument of the exponential function in (\ref{heur:E-conflict}). Therefore, if we choose $\ep$ small enough so that $c(t):=\gamma(m,\ep,\lambda,t)$ is positive, this procedure delivers us a solution  $u(t):=w(m,\ep,\lambda,t)$ to (\ref{pbm:main-ode}) whose energy grows exponentially as the right-hand side of (\ref{heur:E-conflict}). 
\subparagraph{\textmd{\textit{Definition of sequences}}}

Let us choose a sequence $\theta_{n}$ of positive real numbers such that
\begin{equation}
\theta_{n}\to\theta^{-}
\qquad\mbox{as }n\to +\infty.
\nonumber
\end{equation}

Let us consider the sequence $\li$ of the eigenvalues of the operator. Since $2\sigma<\theta<1$, and the sequence $\li$ was assumed to be unbounded, up to passing to a subsequence (not relabeled) we can assume that for every $n\geq 0$ it turns out that
\begin{equation}
\li^{2(1-\theta)}\geq 2,
\label{hp:li-1}
\end{equation}
and 
\begin{equation}
\li^{\theta-2\sigma}\geq 32\delta.
\label{hp:li-2}
\end{equation}

Moreover, for every $n\geq 1$ we can also assume that
\begin{eqnarray}
 & \li^{\theta}\geq 5\lambda_{n-1}^{\theta}, & 
 \label{hp:li-3}   \\[1ex]
 & 2\pi^{\alpha}\li^{2(1-\theta)}\geq 3\lambda_{n-1}^{\theta\alpha}, & 
 \label{hp:li-alpha}   \\[1ex]
 & \pi\li^{\theta}\geq 16\lambda_{n-1}^{\theta}\cdot\left\{(4\delta\li^{2\sigma}+\lambda_{n-1}^{2}+2)\cdot n+3\li^{\theta_{n}}+4(1-\theta)\log\li\right\}. & 
 \label{hp:li-tremenda}
\end{eqnarray}

Indeed, the choice of a (sub)sequence of eigenvalues satisfying these properties can be done inductively. Once that $\lambda_{n-1}$ has been chosen, (\ref{hp:li-3}) and (\ref{hp:li-alpha}) can be easily fulfilled because the sequence of eigenvalues is unbounded. As for (\ref{hp:li-tremenda}), where $\li$ appears on both sides, it is enough to observe that $\theta>2\sigma$, and that $\theta_{n}<\theta$ is a fixed exponent in the moment in which we need to choose $\li$, and therefore the left-hand side grows to $+\infty$ faster than the right-hand side, as a function of $\li$.

Finally, let us set
\begin{equation}
\ep:=\frac{1}{32}, 
\hspace{3em}
m_{n}:=\frac{1}{\li^{1-\theta}},
\hspace{3em}
M_{n}:=m_{n}^{2}+\frac{\delta^{2}}{\lambda_{n}^{2-4\sigma}},
\label{defn:mn}
\end{equation}
and 
\begin{equation}
	t_{n}:=\frac{4\pi}{\li^{\theta}},
	\qquad
	s_{n}:=\frac{\pi}{\li^{\theta}}
	\left\lfloor \frac{2\li^{\theta}}{\liu^{\theta}}\right\rfloor,
	\qquad
	t_{n}':=t_{n}-\frac{\pi}{\li^{\theta}},
	\qquad
	s_{n}':=s_{n}+\frac{\pi}{\li^{\theta}}
	\label{defn:tk-sk}
\end{equation}
where $\lfloor\alpha\rfloor$ denotes the largest integer less than or 
equal to $\alpha$.

\subparagraph{\textmd{\textit{Properties of the sequences}}}

In this paragraph we collect the properties of the sequences that are needed in the sequel.  First of all, the sequence $\li$ is increasing, and $\li\to +\infty$ as $n\to +\infty$. 

From (\ref{defn:mn}) it follows that the sequence $M_{n}$ is decreasing, and it satisfies
\begin{equation}
M_{n}\geq m_{n}^{2}=\frac{1}{\li^{2(1-\theta)}}
\qquad
\forall n\in\n
\label{est:Mi>}
\end{equation}
and, keeping (\ref{hp:li-2}) into account, also
\begin{equation}
M_{n}\leq \frac{3}{2}\,m_{n}^{2}=\frac{3}{2}\,\frac{1}{\li^{2(1-\theta)}}
\qquad
\forall n\in\n.
\label{est:Mi<}
\end{equation}

From (\ref{defn:tk-sk}) and (\ref{hp:li-3}) it follows that
$$s_{n+1}'<t_{n}'<t_{n}<s_{n}<s_{n}'
\quad\quad
\forall n\in\n,$$
and all these sequences tend to 0 as $n\to +\infty$. In addition, for every $n\in\n$ it turns out that
\begin{equation}
	\sin(m_{n}\li t_{n})=\sin(m_{n}\li s_{n})=0,
	\label{est:sin}
\end{equation}
and 
\begin{equation}
	|\cos(m_{n}\li t_{n})|=|\cos(m_{n}\li s_{n})|=1.
	\label{est:cos}
\end{equation}

From (\ref{hp:li-3}) we deduce that for every $n\in\n$ it turns out that
\begin{equation}
s_{n}\geq\frac{\pi}{\lambda_{n-1}^{\theta}},
\hspace{3em}
s_{n}-t_{n}\geq\frac{\pi}{\lambda_{n-1}^{\theta}},
\hspace{3em}
t_{n}'-s_{n+1}'\geq\frac{4\pi}{5\lambda_{n}^{\theta}}.
\label{est:tn-sn}
\end{equation}

Finally, from (\ref{hp:li-tremenda}) and the first inequality in (\ref{est:tn-sn}) we deduce that
\begin{equation}
2\ep\li^{\theta}s_{n}\geq \left(4\delta\li^{2\sigma}+\lambda_{n-1}^{2}+2\right)\cdot n+3\li^{\theta_{n}}+4(1-\theta)\log\li
\qquad
\forall n\geq 1.
\label{est:els-1}
\end{equation}

\subparagraph{\textmd{\textit{Definition of smooth junctions}}}

Plugging $\ep=1/32$ into (\ref{defn:gelt}) we obtain that 
\begin{equation}
\gamma\left(m,\frac{1}{32},\lambda,t\right)=m^{2}+\frac{\delta^{2}}{\lambda^{2-4\sigma}}+m^{2}f(m\lambda t),
\nonumber
\end{equation}
where 
\begin{equation}
f(x):=-\frac{1}{64}\sin^{4}x-\frac{1}{4}\sin(2x)
\nonumber
\end{equation}
is a $\pi$-periodic function such that $f(\pi z)=0$ for every $z\in\z$, and
\begin{equation}
-\frac{1}{2}\leq f\left(x\right)\leq\frac{1}{2}
\qquad
\forall x\in\re.
\label{est:f-bounded}
\end{equation}

In order to create smooth junctions, we choose two function $g_{1}:[0,\pi]\to\re$ and $g_{2}:[0,\pi]\to\re$ of class $C^{\infty}$ such that
\begin{equation}
-\frac{1}{2}\leq g_{1}(x)\leq\frac{1}{2}
\qquad\mbox{and}\qquad
-\frac{1}{2}\leq g_{2}(x)\leq\frac{1}{2}
\label{est:g-bounded}
\end{equation}
for every $x\in[0,\pi]$, and such that the piecewise defined function 
$$\widehat{f}(x):=\left\{
\begin{array}{l@{\qquad}l}
g_{1}(x) & \mbox{if }x\in[0,\pi], \\[0.5ex]
f(x) & \mbox{if }x\in[\pi,2\pi], \\[0.5ex]
g_{2}(x-2\pi) & \mbox{if }x\in[2\pi,3\pi], \\[0.5ex]
0 & \mbox{if }x\not\in(0,3\pi)
\end{array}\right.$$
belongs to $C^{\infty}(\re)$. We observe that, due to the periodicity of $f(x)$, we can repeat the construction above using more blocks of $f(x)$, in the sense that the function
$$\widehat{f}_{j}(x):=\left\{
\begin{array}{l@{\qquad}l}
g_{1}(x) & \mbox{if }x\in[0,\pi], \\[0.5ex]
f(x) & \mbox{if }x\in[\pi,j\pi], \\[0.5ex]
g_{2}(x-j\pi) & \mbox{if }x\in[j\pi,(j+1)\pi], \\[0.5ex]
0 & \mbox{if }x\not\in(0,(j+1)\pi)
\end{array}\right.$$
still belongs to $C^{\infty}(\re)$ for every integer number $j\geq 2$. 

Finally, we choose any function $h:\re\to\re$ of class $C^{\infty}$ such that
\begin{itemize}
  \item $h(x)=0$ for every $x\leq 0$,
  \item $h(x)=1$ for every $x\geq 1$,
  \item $h(x)$ is strictly increasing in $(0,1)$.
\end{itemize}   

\subparagraph{\textmd{\textit{Definition of $c(t)$}}}

Let us define the time-dependent coefficient $c:\re\to\re$. To begin with, for every $n\in\n$ we consider the function $\ell_{n}:\re\to\re$ defined by
\begin{equation}
\ell_{n}(t):= M_{n+1}+(M_{n}-M_{n+1})\cdot h\left(\frac{t-s_{n+1}'}{t_{n}'-s_{n+1}'}\right)
\quad\quad
\forall t\in\re,
\label{defn:li}
\end{equation}
which represents an increasing junction of class $C^{\infty}$ between the constants $M_{n+1}$ and $M_{n}$ in the interval $[s_{n+1}',t_{n}']$. Then in every interval $[s_{n+1}',s_{n}']$ we set
\begin{equation}
c(t):=\left\{
\begin{array}{l@{\qquad}l}
	\ell_{n}(t)   &   \mbox{if }t\in[s_{n+1}',t_{n}'],  \\[1ex]
	M_{n}+m_{n}^{2}g_{1}(m_{n}\li(t-t_{n}'))   &   \mbox{if }t\in[t_{n}',t_{n}],  \\[1ex]
	M_{n}+m_{n}^{2}f(m_{n}\li t) & \mbox{if }t\in[t_{n},s_{n}],   \\[1ex]
	M_{n}+m_{n}^{2}g_{2}(m_{n}\li(t-s_{n}))   &   \mbox{if }t\in[s_{n},s_{n}'].
\end{array}
\right.
\nonumber
\end{equation}

\begin{figure}[htbp]
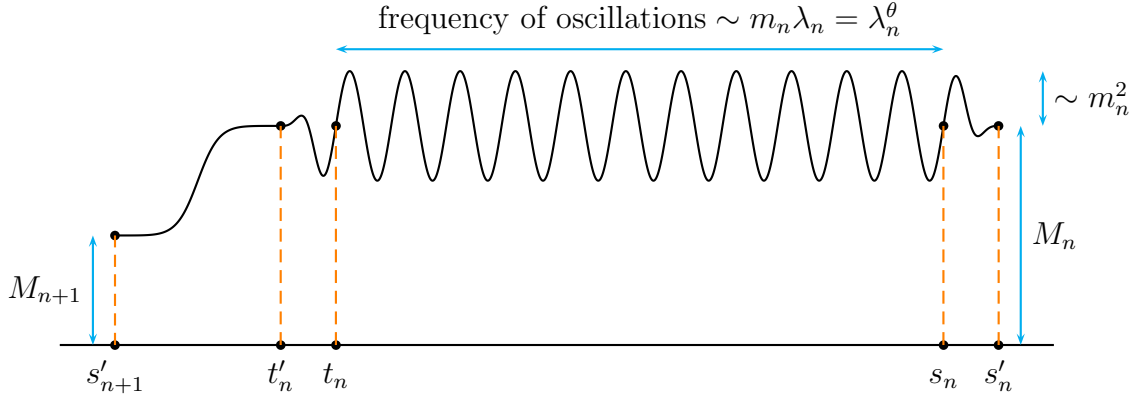

\begin{center}
\psset{unit=4ex}
\pspicture(-2,-2)(17,5.5)

\psplot[plotpoints=180]{0}{3}{x 4 exp x 4 exp x 3 sub 4 exp add div 2 mul 1 add}
\psplot[plotpoints=180]{3}{4}{x 3 sub 2 exp 360 x mul sin mul x 3 sub 2 exp x 4 sub 2 exp add div 3 add}
\psplot[plotpoints=1800]{4}{15}{360 x mul sin 3 add}
\psplot[plotpoints=180]{15}{16}{x 16 sub 2 exp 360 x mul sin mul x 16 sub 2 exp x 15 sub 2 exp add div 3 add}

\psline(-1,-1)(17,-1)
\psdots(0,-1)(3,-1)(4,-1)(15,-1)(16,-1)
\psdots(0,1)(3,3)(4,3)(15,3)(16,3)

\rput[B](0,-1.7){$s_{n+1}'$}
\rput[B](3,-1.7){$t_{n}'$}
\rput[B](4,-1.7){$t_{n}$}
\rput[B](15,-1.7){$s_{n}$}
\rput[B](16,-1.7){$s_{n}'$}

\psline[linecolor=orange,linestyle=dashed](0,-1)(0,1)
\psline[linecolor=orange,linestyle=dashed](3,-1)(3,3)
\psline[linecolor=orange,linestyle=dashed](4,-1)(4,3)
\psline[linecolor=orange,linestyle=dashed](15,-1)(15,3)
\psline[linecolor=orange,linestyle=dashed](16,-1)(16,3)

\pcline[linecolor=cyan]{<->}(16.4,-1)(16.4,3)
\rput[l](16.6,1){$M_{n}$}
\pcline[linecolor=cyan]{<->}(-0.4,-1)(-0.4,1)
\rput[r](-0.6,0){$M_{n+1}$}
\pcline[linecolor=cyan]{<->}(16.8,3)(16.8,4)
\rput[l](17,3.5){$\sim m_{n}^{2}$}
\pcline[linecolor=cyan]{<->}(4,4.4)(15,4.4)
\lput{U}{\uput[90](0,0){{$\mbox{frequency of oscillations}\sim m_{n}\lambda_{n}=\lambda_{n}^{\theta}$}}}

\endpspicture
\caption{basic block of $c(t)$ between $s_{n+1}'$ and $s_{n}'$}
\label{fig:c(t)-1}
\end{center}
\end{figure}

Figure~\ref{fig:c(t)-1} describes the shape of $c(t)$ in the interval $[s_{n+1}',s_{n}']$. This is the building block of the entire construction. We observe in particular that\begin{equation}
c(t_{n}')=c(t_{n})=c(s_{n})=c(s_{n}')=M_{n}
\qquad
\forall n\in\n.
\label{eqn:c-endpoints}
\end{equation}

\begin{figure}[htbp]
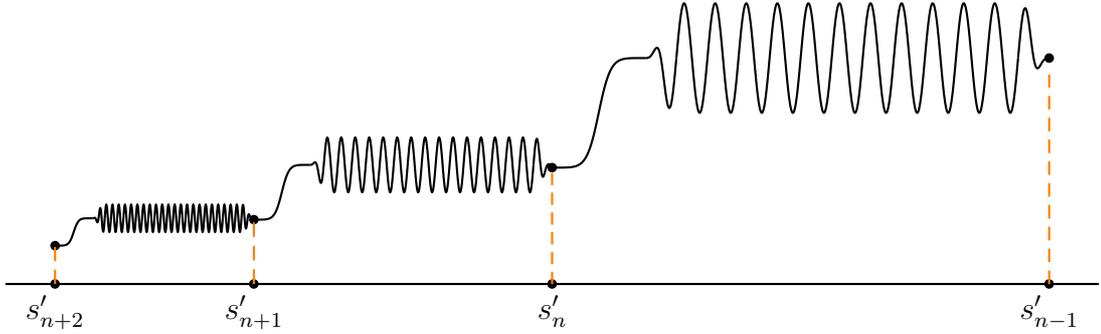

\begin{center}
\psset{unit=3.6ex}
\pspicture(-2,-1.5)(21,5.5)

\rput(0,0){
\psset{xunit=0.9ex,yunit=1ex}
\psplot[plotpoints=180]{0}{3}{x 4 exp x 4 exp x 3 sub 4 exp add div 2 mul 1 add}
\psplot[plotpoints=180]{3}{4}{x 3 sub 2 exp 720 x mul sin mul x 3 sub 2 exp x 4 sub 2 exp add div 3 add}
\psplot[plotpoints=1000]{4}{15}{720 x mul sin 3 add}
\psplot[plotpoints=180]{15}{16}{x 16 sub 2 exp 720 x mul sin mul x 16 sub 2 exp x 15 sub 2 exp add div 3 add}
\psdot(0,1)
}

\rput(4,0.25){
\psset{xunit=1.35ex,yunit=2ex}
\psplot[plotpoints=180]{0}{3}{x 4 exp x 4 exp x 3 sub 4 exp add div 2 mul 1 add}
\psplot[plotpoints=180]{3}{4}{x 3 sub 2 exp 480 x mul sin mul x 3 sub 2 exp x 4 sub 2 exp add div 3 add}
\psplot[plotpoints=1000]{4}{15}{480 x mul sin 3 add}
\psplot[plotpoints=180]{15}{16}{x 16 sub 2 exp 480 x mul sin mul x 16 sub 2 exp x 15 sub 2 exp add div 3 add}
\psdot(0,1)
}

\rput(10,0.75){
\psset{xunit=2.25ex,yunit=4ex}
\psplot[plotpoints=180]{0}{3}{x 4 exp x 4 exp x 3 sub 4 exp add div 2 mul 1 add}
\psplot[plotpoints=180]{3}{4}{x 3 sub 2 exp 360 x mul sin mul x 3 sub 2 exp x 4 sub 2 exp add div 3 add}
\psplot[plotpoints=1000]{4}{15}{360 x mul sin 3 add}
\psplot[plotpoints=180]{15}{16}{x 16 sub 2 exp 360 x mul sin mul x 16 sub 2 exp x 15 sub 2 exp add div 3 add}
\psdot(0,1)
\psdot(16,3)
}

\psline(-1,-0.5)(21,-0.5)
\psdots(0,-0.5)(4,-0.5)(10,-0.5)(20,-0.5)
\rput[B](0,-1.2){$s_{n+2}'$}
\rput[B](4,-1.2){$s_{n+1}'$}
\rput[B](10,-1.2){$s_{n}'$}
\rput[B](20,-1.2){$s_{n-1}'$}

\psline[linestyle=dashed,linecolor=orange](0,-0.5)(0,0.25)
\psline[linestyle=dashed,linecolor=orange](4,-0.5)(4,0.75)
\psline[linestyle=dashed,linecolor=orange](10,-0.5)(10,1.75)
\psline[linestyle=dashed,linecolor=orange](20,-0.5)(20,3.75)

\endpspicture
\caption{glueing of blocks in the definition of $c(t)$}
\label{fig:c(t)-2}
\end{center}
\end{figure}

The building blocks are repeated as described in Figure~\ref{fig:c(t)-2}. This defines $c(t)$ in the interval $(0,s_{0}']$. We complete the definition by setting $c(t):=0$ for every $t\leq 0$, and $c(t):= M_{0}$ for every $t\geq s_{0}'$. The resulting function $c(t)$ belongs to $C^{\infty}(\re\setminus\{0\})$, and its derivatives of any order vanish in the points $s_{n}'$.

In words, the idea of the construction is the following.
\begin{itemize}

\item In the interval $[t_{n},s_{n}]$ the function $c(t)$ coincides with $\gamma(m_{n},\ep,\li,t)$. This function oscillates, with frequency of order $\lambda_{n}^{\theta}$ and amplitude of order $m_{n}^{2}$, around the mean value $M_{n}$. Both the mean value $M_{n}$ and the amplitude of oscillations tend to~0 as $n\to +\infty$, while the frequency of oscillations diverges to $+\infty$.

\item In the intervals $[t_{n}',t_{n}]$ and $[s_{n},s_{n}']$ the function $c(t)$ is a $C^{\infty}$ junction between the oscillating function of the interval $[t_{n},s_{n}]$ and the constant $M_{n}$. 

\item  In the interval $[s_{n+1}',t_{n}']$ the function $c(t)$ is an increasing junction of class $C^{\infty}$ between the constant $M_{n+1}$ and the constant $M_{n}$.

\end{itemize}

We point out that the key feature of the construction is the highly oscillatory behavior of $c(t)$ in the intervals $[t_{n},s_{n}]$; all the rest is aimed at creating a smooth transition between the values of $c(t)$ in these intervals.

\subparagraph{\textmd{\textit{Definition of $u(t)$}}}

For every $n\in\n$, we consider the solution $u_{n}(t)$ to 
the ordinary differential equation
\begin{equation}
u_{n}''(t)+2\delta\li^{2\sigma}u_{n}'(t)+\li^{2}c(t)u_{n}(t)=0,
\label{uk-eqn}
\end{equation}
with ``initial'' data
\begin{equation}
	u_{n}(t_{n})=0,
	\quad\quad
	u_{n}'(t_{n})=
	m_{n}\li\exp\left((2\ep m_{n}\li-\delta\li^{2\sigma})t_{n}\right).
	\label{uk-data}
\end{equation}

Then we set
\begin{equation}
	a_{n}:=\frac{m_{n}}{(n+1)\li^{\theta}}\exp(-\li^{\theta_{n}}),
	\label{defn:ak}
\end{equation}
and we consider the solution $u(t)$ to (\ref{pbm:eqn}) defined by
$$u(t):=\sum_{n=0}^{\infty}a_{n}u_{n}(t)e_{n}.$$

We claim that $c(t)$ satisfies (\ref{th:c-k-a}) and (\ref{th:c-dh}), and that $u(t)$ satisfies (\ref{th:u0-reg}) and (\ref{th:ut-schifo}). The rest of the proof is a verification of these claims.

\subparagraph{\textmd{\textit{Continuity and degenerate hyperbolicity of $c(t)$}}}

We prove that for every $n\in\n$ it turns out that
\begin{equation}
\frac{1}{2\li^{2(1-\theta)}}\leq c(t)\leq\frac{2}{\li^{2(1-\theta)}}
\qquad
\forall t\in[t_{n}',s_{n}'],
\label{est:c-1}
\end{equation}
and
\begin{equation}
\frac{1}{\lambda_{n+1}^{2(1-\theta)}}\leq c(t)\leq\frac{2}{\li^{2(1-\theta)}}
\qquad
\forall t\in[s_{n+1}',t_{n}'].
\label{est:c-2}
\end{equation}

If we prove these estimates, then from (\ref{hp:li-1}) it follows that $c(t)$ satisfies the degenerate hyperbolicity assumption (\ref{th:c-dh}). Moreover, since $\li\to+\infty$, we obtain also that $c(t)\to 0$ as $t\to 0^{+}$, which proves the continuity of $c(t)$ on the whole real line.

In order to prove (\ref{est:c-1}), from (\ref{est:f-bounded}) and (\ref{est:g-bounded}) we obtain the estimates
\begin{equation}
M_{n}-\frac{1}{2}m_{n}^{2}\leq c(t)\leq M_{n}+\frac{1}{2}m_{n}^{2}
\qquad
\forall t\in[t_{n}',s_{n}'],
\nonumber
\end{equation}
so that the conclusion follows from (\ref{est:Mi>}) and (\ref{est:Mi<}).

In order to prove (\ref{est:c-2}), we just recall that in the interval $[s_{n+1}',t_{n}']$ the function $c(t)$ is a smooth increasing junction between the constants $M_{n+1}$ and $M_{n}$, and we conclude by exploiting again (\ref{est:Mi>}) and (\ref{est:Mi<}).

\subparagraph{\textmd{\textit{Regularity of $c(t)$ and estimates for its derivatives}}}

For every positive integer $j$, let $c^{(j)}(t)$ denote the $j$-th derivative of $c(t)$. To begin with, we prove that there exists a constant $\Gamma_{j}$ such that
\begin{equation}
|c^{(j)}(t)|\leq\Gamma_{j}\li^{(2+j)\theta-2}
\qquad
\forall n\in\n,\quad\forall t\in[s_{n+1}',s_{n}'].
\label{est:cj}
\end{equation}

The constant $\Gamma_{j}$ depends only on the $L^{\infty}$ norms of the derivatives of order $j$ of the functions $f(x)$, $g_{1}(x)$, $g_{2}(x)$, and $h(x)$. 

If $k\geq 1$, then $(2+j)\theta<2$ for every $j\leq k$, and hence (\ref{est:cj}) implies that $c^{(j)}(t)\to 0$ as $t\to 0^{+}$ for every $j\leq k$, which proves that $c\in C^{k}(\re)$.

In order to establish (\ref{est:cj}), it is enough to examine the definition of $c(t)$ in the four subintervals whose union is $[s_{n+1}',s_{n}']$.
\begin{itemize}

\item  In the interval $[t_{n},s_{n}]$ the function $c(t)$ is obtained from $f(x)$ through horizontal and vertical rescaling, and therefore
$$|c^{(j)}(t)|=m_{n}^{2}(m_{n}\li)^{j}\left|f^{(j)}(m_{n}\li t)\right|\leq\li^{(2+j)\theta-2}\cdot\|f^{(j)}\|_{\infty}.$$

\item  In the intervals $[t_{n}',t_{n}]$ and $[s_{n},s_{n}']$ the argument is exactly the same, just with $g_{1}$ and $g_{2}$ instead of $f$.

\item  In the interval $[s_{n+1}',t_{n}']$ the function $c(t)$ is a rescaling of the function $h(x)$. Thus from (\ref{defn:li}) it follows that
$$|c^{(j)}(t)|=\frac{M_{n}-M_{n+1}}{(t_{n}'-s_{n+1}')^{j}}\cdot\left|h^{(j)}\left(\frac{t-s_{n+1}'}{t_{n}'-s_{n+1}'}\right)\right|\leq\frac{M_{n}-M_{n+1}}{(t_{n}'-s_{n+1}')^{j}}\cdot\|h^{(j)}\|_{\infty}.$$

Now we estimate the numerator with $M_{n}$, which in turn we estimate as in (\ref{est:Mi<}), and we estimate the denominator as in third inequality in (\ref{est:tn-sn}). This is enough to conclude (\ref{est:cj}) also in this case.

\end{itemize}

Now we show that actually $c\in C^{k,\alpha}(\re)$. Since $c^{(k)}(t)$ is continuous, and constant for $t\leq 0$ and $t\geq s_{0}'$, it is enough to prove the $\alpha$-H\"older continuity of $c^{(k)}(t)$ in $(0,s_{0}')$. 

To this end, we can limit ourselves to showing that the $\alpha$-H\"older constant of $c^{(k)}(t)$ in the interval $[s_{n+1}',s_{n}']$ is bounded from above independently of $n$, and 
\begin{equation}
\left|c^{(k)}(s_{i}')-c^{(k)}(s_{j}')\right|\leq\left|s_{i}'-s_{j}'\right|^{\alpha}
\qquad
\forall (i,j)\in\n^{2}.
\label{ck-sj}
\end{equation}

Indeed, let us consider any interval $[x,y]\subseteq(0,s_{0}')$. If $x$ and $y$ lie in the same interval of the form $[s_{n+1}',s_{n}']$, then $|c^{(k)}(y)-c^{(k)}(x)|$ can be controlled in terms of $|y-x|^{\alpha}$ because of the uniform bound on the H\"older constants. The same is true if $x$ and $y$ lie in neighboring intervals. In the remaining case, there exist two positive indices $i<j$ such that
\begin{equation}
s_{j+1}'<x\leq s_{j}'<s_{i}'\leq y<s_{i-1}'.
\nonumber
\end{equation}

In this case we write
$$|c^{(k)}(y)-c^{(k)}(x)|\leq|c^{(k)}(y)-c^{(k)}(s_{i}')|+|c^{(k)}(s_{i}')-c^{(k)}(s_{j}')|+|c^{(k)}(s_{j}')-c^{(k)}(x)|,$$
and we observe that the central term is less that $|y-x|^{\alpha}$ because of (\ref{ck-sj}), while the other two terms can be controlled by exploiting once again the uniform bound on the H\"older constants.

\begin{itemize}

\item As for the uniform estimate of the H\"older constant, the same scaling arguments used in the estimates of the derivatives of $c(t)$ shows that this constant is less than or equal to
$$\Gamma_{k,\alpha}\,\li^{(2+k+\alpha)\theta-2},$$
where $\Gamma_{k,\alpha}$ is proportional to the $\alpha$-H\"older constants of the $k$-th derivatives of the functions $f(x)$, $g_{1}(x)$, $g_{2}(x)$, and $h(x)$.  Due to (\ref{defn:theta}), the exponent of $\li$ is zero, and hence the bound is independent of $n$. 

\item  As for (\ref{ck-sj}), we observe that the inequality is trivial if $k\geq 1$ because all derivatives of $c(t)$ vanish in the points $s_{n}'$. If $k=0$, we assume without loss of generality that $i<j$, and from (\ref{eqn:c-endpoints}) and (\ref{est:Mi<}) we deduce that
$$|c(s_{i}')-c(s_{j}')|=M_{i}-M_{j}\leq M_{i}\leq\frac{3}{2}\,\frac{1}{\lambda_{i}^{2(1-\theta)}},$$
while from the second inequality in (\ref{est:tn-sn}) we deduce that
$$s_{i}'-s_{j}'\geq s_{i}-t_{i}\geq\frac{\pi}{\lambda_{i-1}^{\theta}}.$$

At this point (\ref{ck-sj}) follows from (\ref{hp:li-alpha}).

\end{itemize} 

\subparagraph{\textmd{\textit{Energy functions}}}

Let us consider the classic energy functions
$$E_{n}(t):=|u_{n}'(t)|^{2}+m_{n}^{2}\li^{2}|u_{n}(t)|^{2},$$
$$F_{n}(t):=|u_{n}'(t)|^{2}+\li^{2}c(t)|u_{n}(t)|^{2}.$$

Since $m_{n}\leq 1$, and $0<c(t)\leq 1$ for every $t\geq 0$, it turns out that
\begin{equation}
|u_{n}'(t)|^{2}+\li^{2}|u_{n}(t)|^{2}\leq\frac{1}{m_{n}^{2}}E_{n}(t)
\qquad
\forall n\in\n,\quad\forall t\geq 0,
\nonumber
\end{equation}
\begin{equation}
|u_{n}'(t)|^{2}+\li^{2}|u_{n}(t)|^{2}\geq F_{n}(t)
\qquad
\forall n\in\n,\quad\forall t\geq 0.
\nonumber
\end{equation}

Therefore, (\ref{th:u0-reg}) is proved if we show that
\begin{equation}
	\sum_{n=0}^{\infty}
	a_{n}^{2}\frac{1}{m_{n}^{2}}E_{n}(0)\exp\left(2r\li^{1/s}\right)<+\infty
	\quad\quad
	\forall r>0,\quad\forall s>1+\frac{k+\alpha}{2},
	\label{th:u0-equiv}
\end{equation}
while (\ref{th:ut-schifo}) is proved if we show that for every $t>0$ it turns out that
\begin{equation}
	\sum_{n=0}^{\infty}
	a_{n}^{2}F_{n}(t)\exp\left(-2R\li^{1/S}\right)=+\infty
	\quad\quad
	\forall R>0,\quad\forall S>1+\frac{k+\alpha}{2}.
	\label{th:ut-equiv}
\end{equation}

We can neglect the Sobolev parameter $\beta$ in the spaces involved in (\ref{th:u0-reg}) and (\ref{th:ut-schifo}) because powers of $\li$ are lower order terms with respect to the exponentials. Thus in the sequel we just have to estimate $E_{n}(0)$ and $F_{n}(t)$.

\subparagraph{\textmd{\textit{Energy estimates in $[0,t_{n}]$}}}

We prove that 
\begin{equation}
E_{n}(0)\leq\li^{2\theta}\exp(5\pi)
\qquad
\forall n\in\n.
\label{th:Ei(0)}
\end{equation}

To begin with, from (\ref{uk-data}), (\ref{defn:mn}) and (\ref{defn:tk-sk}) we obtain that
\begin{equation}
E_{n}(t_{n})=|u_{n}'(t_{n})|^{2}\leq m_{n}^{2}\li^{2}\exp(4m_{n}\li\ep t_{n})=\li^{2\theta}\exp(\pi/2).
\label{est:Ei(ti)}
\end{equation}

Moreover, the time-derivative of $E_{n}(t)$ can be estimated as
\begin{eqnarray}
E_{n}'(t)  &  =  &  -4\delta\li^{2\sigma}|u_{n}'(t)|^{2}-m_{n}\li\left(\frac{c(t)}{m_{n}^{2}}-1\right)\cdot 2m_{n}\li u_{n}(t)u_{n}'(t)
\nonumber
\\[1ex]
  &  \geq  &  -4\delta\li^{2\sigma}E_{n}(t)-m_{n}\li\left|\frac{c(t)}{m_{n}^{2}}-1\right|E_{n}(t).
  \label{est:En'}
\end{eqnarray}

Since $\li$ is increasing, from (\ref{est:c-1}) and (\ref{est:c-2}) we obtain that
$$0\leq c(t)\leq 2m_{n}^{2}
\qquad
\forall t\in[0,t_{n}],$$
and therefore from (\ref{est:En'}) we deduce that
$$E_{n}'(t)\geq -\left(4\delta\li^{2\sigma}+\li^{\theta}\right)E_{n}(t)
\qquad
\forall t\in[0,t_{n}].$$

Integrating this differential inequality, and keeping (\ref{defn:tk-sk}) and (\ref{hp:li-2}) into account, we deduce that
\begin{equation}
E_{n}(0)\leq E_{n}(t_{n})\exp\left(4\delta\li^{2\sigma}t_{n}+\li^{\theta} t_{n}\right)\leq E_{n}(t_{n})\exp(\pi/2+4\pi).
\label{est:Ei(0)-Ei(ti)}
\end{equation}

Plugging (\ref{est:Ei(ti)}) into (\ref{est:Ei(0)-Ei(ti)}), we obtain (\ref{th:Ei(0)}).

\subparagraph{\textmd{\textit{Energy estimates in $[t_{n},s_{n}]$}}}

In this interval the solution to (\ref{uk-eqn})--(\ref{uk-data}) is given by the explicit formula $u_{n}(t):=w(m_{n},\ep,\li,t)$, where
$w(m,\ep,\lambda,t)$ is the function defined in (\ref{defn:welt}). Keeping (\ref{est:sin}) and (\ref{est:cos}) into account, we deduce that $u_{n}(s_{n})=0$ and
$$|u_{n}'(s_{n})|=m_{n}\li\exp\left((2m_{n}\li\ep-\delta\li^{2\sigma})s_{n}\right)=\li^{\theta}\exp\left((2\li^{\theta}\ep-\delta\li^{2\sigma})s_{n}\right).$$

Therefore, from (\ref{hp:li-2}) and the definition of $\ep$ it follows that
$$|u_{n}'(s_{n})|\geq\li^{\theta}\exp\left(\ep \li^{\theta}s_{n}\right),$$
and hence
\begin{equation}
	F_{n}(s_{n})=E_{n}(s_{n})=|u_{n}'(s_{n})|^{2}\geq \li^{2\theta}\exp(2\ep\li^{\theta}s_{n})
	\qquad
	\forall n\in\n.
	\label{est:fk-sk}
\end{equation}

\subparagraph{\textmd{\textit{Energy estimates in $[s_{n},+\infty)$}}}

We prove that for every $t\geq s_{n}$ it turns out that
\begin{eqnarray}
F_{n}(t) & \geq & \li^{2\theta}\exp\left(2\ep\li^{\theta}s_{n}\right)\cdot 
\nonumber \\[1ex]
 & & \mbox{}\cdot\exp\left(-(4\delta\li^{2\sigma}+2\Gamma_{1}\lambda_{n-1}^{2})t-2\pi\Gamma_{1}-2(1-\theta)\log\li\strut\right),
\label{est:fk-t}
\end{eqnarray}
where $\Gamma_{1}$ is the constant for which (\ref{est:cj}) holds true in the case $j=1$.

To begin with, we estimate the time-derivative of the hyperbolic energy as
\begin{eqnarray*}
F_{n}'(t) & = & -4\delta\li^{2\sigma}|u_{n}'(t)|^{2}+\li^{2}c'(t)|u_{n}(t)|^{2}\\[1ex]
 & \geq & -4\delta\li^{2\sigma}|u_{n}'(t)|^{2}-\frac{|c'(t)|}{c(t)}\cdot \li^{2}c(t)|u_{n}(t)|^{2}  \\[1ex]
 & \geq & -\left(4\delta\li^{2\sigma}+\frac{|c'(t)|}{c(t)}\right)F_{n}(t).
\end{eqnarray*}

Integrating this differential inequality, we obtain that
\begin{equation}
F_{n}(t)\geq F_{n}(s_{n})\exp\left(-4\delta\li^{2\sigma}t-\int_{s_{n}}^{t}\frac{|c'(s)|}{c(s)}\,ds\right)
\qquad
\forall t\geq s_{n}.
\label{est:Fi-int}
\end{equation}

In order to estimate the last integral, we write it as the sum of three terms
$$I_{1}:=\int_{s_{n}}^{s_{n}'}\frac{|c'(s)|}{c(s)}\,ds,
\hspace{3em}
I_{2}:=\int_{s_{n}'}^{t_{n-1}'}\frac{|c'(s)|}{c(s)}\,ds,
\hspace{3em}
I_{3}:=\int_{t_{n-1}'}^{t}\frac{|c'(s)|}{c(s)}\,ds,$$
which we consider separately (we assume to be in the worst case scenario where $t>t_{n-1}'$, so that all the integrals need to be estimated).

\begin{itemize}

\item In the interval $[s_{n},s_{n}']$ we deduce from (\ref{est:c-1}) and (\ref{est:cj}) that
$$c(t)\geq\frac{1}{2\li^{2(1-\theta)}}
\qquad\mbox{and}\qquad
|c'(t)|\leq\Gamma_{1}\li^{3\theta-2},$$
from which we conclude that
\begin{equation}
I_{1}\leq \Gamma_{1}\li^{3\theta-2}\cdot2\li^{2(1-\theta)}\cdot(s_{n}'-s_{n})\leq 2\pi\Gamma_{1}.
\label{est:I1}
\end{equation}

\item  In the interval $[s_{n}',t_{n-1}']$ the function $c(t)$ is an increasing junction between $M_{n}$ and $M_{n-1}$, and therefore
$$\int_{s_{n}'}^{t_{n-1}'}\frac{|c'(s)|}{c(s)}\,ds=\int_{s_{n}'}^{t_{n-1}'}\frac{c'(s)}{c(s)}\,ds=\log\frac{c(t_{n-1}')}{c(s_{n}')}=\log\frac{M_{n-1}}{M_{n}}.$$

From (\ref{hp:li-1}), (\ref{est:Mi>}) and (\ref{est:Mi<}) we know that $M_{n-1}\leq 1$ and $M_{n}\geq\li^{-2(1-\theta)}$, and hence
\begin{equation}
I_{2}=\log\frac{M_{n-1}}{M_{n}}\leq\log(\li^{2(1-\theta)})=2(1-\theta)\log\li.
\label{est:I2}
\end{equation}

\item  Let us consider the interval $[t_{n-1}',t]$, and let us observe that
$$[t_{n-1}',t]\subseteq[t_{n-1}',s_{n-1}']\cup\bigcup_{i=0}^{n-2}[s_{i+1}',s_{i}']\cup[s_{0}',+\infty).$$

Due to the estimates from below in (\ref{est:c-1}) and (\ref{est:c-2}), and recalling that the sequence $\li$ is increasing, we deduce that
\begin{equation}
c(s)\geq\frac{1}{2\lambda_{n-1}^{2(1-\theta)}}
\qquad
\forall s\geq t_{n-1}'.
\nonumber
\end{equation}

As for $|c'(s)|$, we exploit (\ref{est:cj}) with $j=1$, and we obtain that
\begin{equation}
|c'(s)|\leq\max_{i\leq n-1}\Gamma_{1}\lambda_{i}^{3\theta-2}
\qquad
\forall s\geq t_{n-1}'.
\nonumber
\end{equation}

The value of the maximum depends on $\theta$.
\begin{itemize}

\item  If $\theta<2/3$, which corresponds to the case $k\geq 1$, the maximum is attained for $i=0$. Since $\lambda_{0}\geq 1$, the maximum can be estimated from above with $\Gamma_{1}$, and hence
\begin{equation}
\frac{|c'(s)|}{c(s)}\leq 2\Gamma_{1}\lambda_{n-1}^{2(1-\theta)}.
\nonumber
\end{equation}

\item  If $\theta\geq 2/3$, which corresponds to the case $k=0$, the maximum is attained for $i=n-1$, and hence
\begin{equation}
\frac{|c'(s)|}{c(s)}\leq 2\Gamma_{1}\lambda_{n-1}^{\theta}.
\nonumber
\end{equation}

\end{itemize}

In conclusion, in both cases we have proved that
$$\frac{|c'(s)|}{c(s)}\leq 2\Gamma_{1}\lambda_{n-1}^{\max\{2(1-\theta),\theta\}}\leq 2\Gamma_{1}\lambda_{n-1}^{2},$$
and hence
\begin{equation}
I_{3}=\int_{t_{n-1}'}^{t}\frac{|c'(s)|}{c(s)}\,ds\leq 2\Gamma_{1}\lambda_{n-1}^{2}t.
\label{est:I3}
\end{equation}

\end{itemize}

Plugging (\ref{est:I1}), (\ref{est:I2}) and (\ref{est:I3}) into (\ref{est:Fi-int}) we conclude that
$$F_{n}(t)\geq F_{n}(s_{n})\exp\left(-4\delta\li^{2\sigma}t-2\pi\Gamma_{1}-2(1-\theta)\log\li-2\Gamma_{1}\lambda_{n-1}^{2}t\right).$$

Keeping (\ref{est:fk-sk}) into account, we obtain exactly (\ref{est:fk-t}).

\subparagraph{\textmd{\textit{Conclusion}}}

We are now ready to verify (\ref{th:u0-equiv}) and
(\ref{th:ut-equiv}).  Indeed from (\ref{defn:ak}) and (\ref{th:Ei(0)})
it turns out that
\begin{eqnarray*}
	\frac{a_{n}^{2}}{m_{n}^{2}}\cdot E_{n}(0)\cdot\exp(2r\li^{1/s}) & \leq & 
	\frac{1}{(n+1)^{2}\li^{2\theta}}\exp(-2\li^{\theta_{n}})\cdot
	\li^{2\theta}\exp(5\pi)\cdot\exp(2r\li^{1/s})    \\[1ex]
	 & = & \frac{1}{(n+1)^{2}}\exp\left(5\pi-2\li^{\theta_{n}}+2r\li^{1/s}\right).
\end{eqnarray*}

Since $\theta_{n}\to\theta>1/s$, the argument of the exponential is bounded from above, and hence the series in (\ref{th:u0-equiv}) converges.

Let us consider now (\ref{th:ut-equiv}). From (\ref{est:fk-t}) and (\ref{est:els-1}) it follows that
$$F_{n}(t)\geq\li^{2\theta}\exp\left(3\li^{\theta_{n}}+n+2(1-\theta)\log\li\right)=\frac{\li^{2\theta}}{m_{n}^{2}}\exp\left(3\li^{\theta_{n}}+n\right)$$
when $n$ is large enough, 
and therefore
\begin{eqnarray*}
a_{n}^{2}\cdot F_{n}(t)\cdot\exp\left(-2R\li^{1/S}\right) & \geq & \frac{m_{n}^{2}}{(n+1)^{2}\li^{2\theta}}\exp\left(-2\li^{\theta_{n}}\right)\cdot \\[1ex]
& & \mbox{}\cdot\frac{\li^{2\theta}}{m_{n}^{2}}\exp\left(3\li^{\theta_{n}}+n\right)\cdot\exp\left(-2R\li^{1/S}  \right)  \\[1ex]
 & = & \frac{1}{(n+1)^{2}}\exp\left(\li^{\theta_{n}}-2R\li^{1/S}+n\right)
\end{eqnarray*}
for the same values of $n$. Since $\theta_{n}\to\theta>1/S$, the argument of the exponential is eventually greater than $n$, and therefore the series in (\ref{th:ut-equiv}) diverges.\qed



\begin{thebibliography}{10}
\providecommand{\url}[1]{\texttt{#1}}
\providecommand{\urlprefix}{URL }
\providecommand{\selectlanguage}[1]{\relax}
\providecommand{\eprint}[2][]{\url{#2}}

\bibitem{Bui-Reissig}
\textsc{T.~B.~N. Bui}, \textsc{M.~Reissig}.
\newblock The interplay between time-dependent speed of propagation and
  dissipation in wave models.
\newblock In \emph{Fourier analysis}, Trends Math., pages 9--45.
  Birkh\"auser/Springer, Cham, 2014.

\bibitem{dgcs}
\textsc{F.~Colombini}, \textsc{E.~De~Giorgi}, \textsc{S.~Spagnolo}.
\newblock Sur les \'equations hyperboliques avec des coefficients qui ne
  d\'ependent que du temps.
\newblock \emph{Ann. Scuola Norm. Sup. Pisa Cl. Sci. (4)} \textbf{6} (1979),
  no.~3, 511--559.

\bibitem{cjs}
\textsc{F.~Colombini}, \textsc{E.~Jannelli}, \textsc{S.~Spagnolo}.
\newblock Well-posedness in the {G}evrey classes of the {C}auchy problem for a
  nonstrictly hyperbolic equation with coefficients depending on time.
\newblock \emph{Ann. Scuola Norm. Sup. Pisa Cl. Sci. (4)} \textbf{10} (1983),
  no.~2, 291--312.

\bibitem{DAbbicco-Ebert}
\textsc{M.~D'Abbicco}, \textsc{M.~R. Ebert}.
\newblock A class of dissipative wave equations with time-dependent speed and
  damping.
\newblock \emph{J. Math. Anal. Appl.} \textbf{399} (2013), no.~1, 315--332.

\bibitem{Ebert-Reissig}
\textsc{M.~R. Ebert}, \textsc{M.~Reissig}.
\newblock Theory of damped wave models with integrable and decaying in time
  speed of propagation.
\newblock \emph{J. Hyperbolic Differ. Equ.} \textbf{13} (2016), no.~2,
  417--439.

\bibitem{gg:roots}
\textsc{M.~Ghisi}, \textsc{M.~Gobbino}.
\newblock Higher order {G}laeser inequalities and optimal regularity of roots
  of real functions.
\newblock \emph{Ann. Sc. Norm. Super. Pisa Cl. Sci. (5)} \textbf{12} (2013),
  no.~4, 1001--1021.

\bibitem{gg:dgcs-strong}
\textsc{M.~Ghisi}, \textsc{M.~Gobbino}.
\newblock Linear wave equations with time-dependent propagation speed and
  strong damping.
\newblock \emph{J. Differential Equations} \textbf{260} (2016), no.~2,
  1585--1621.

\bibitem{ggh:strong-damping}
\textsc{M.~Ghisi}, \textsc{M.~Gobbino}, \textsc{A.~Haraux}.
\newblock Local and global smoothing effects for some linear hyperbolic
  equations with a strong dissipation.
\newblock \emph{Trans. Amer. Math. Soc.} \textbf{368} (2016), no.~3,
  2039--2079.

\bibitem{Hirosawa-Bui}
\textsc{F.~Hirosawa}, \textsc{T.~B.~N. Bui}.
\newblock On energy estimates for second order hyperbolic equations with {L}evi
  conditions for higher order regularity.
\newblock \emph{Ann. Univ. Ferrara Sez. VII Sci. Mat.} \textbf{57} (2011),
  no.~2, 317--339.

\bibitem{reed}
\textsc{M.~Reed}, \textsc{B.~Simon}.
\newblock \emph{Methods of modern mathematical physics. {I}}.
\newblock Academic Press, Inc. [Harcourt Brace Jovanovich, Publishers], New
  York, second edition, 1980.
\newblock Functional analysis.

\end{thebibliography}

\label{NumeroPagine}

\end{document}